\documentclass[USenglish]{article}

\usepackage[USenglish]{babel}
\usepackage{color}
\usepackage{amsmath,amsthm,amssymb,amsfonts}
\usepackage[margin=2.5 cm]{geometry}
\usepackage[unicode,colorlinks=true,linkcolor=blue,citecolor=green,urlcolor=black,pagebackref=false]{hyperref}
\usepackage{stmaryrd}
\usepackage{tikz}
\usepackage{pgfplots}
\usepackage{mathabx}
\usepackage{custom_macros_general, custom_macros_spaceops, custom_macros_funcs, custom_macros_freqops}

\let\oldmu\mu
\let\oldeps\varepsilon
\let\oldnu\nu
\let\oldlambda\lambda
\let\oldzeta\zeta

\renewcommand{\mu}{\boldsymbol{\oldmu}}
\renewcommand{\varepsilon}{\boldsymbol{\oldeps}}
\renewcommand{\nu}{\boldsymbol{\oldnu}}
\renewcommand{\lambda}{\boldsymbol{\oldlambda}}
\renewcommand{\zeta}{\boldsymbol{\oldzeta}}

\def\revision#1{{#1}}

\def\revtwo#1{{#1}}
\numberwithin{equation}{section}

\DeclareMathAlphabet\mathbfcal{OMS}{cmsy}{b}{n}
\makeatletter
\let\@fnsymbol\@arabic
\makeatother

\title{Wavenumber-explicit $hp$-FEM analysis of Maxwell's equations with impedance boundary conditions in piecewise smooth media}

\author{Jens Markus Melenk\thanks{(melenk@tuwien.ac.at), Institut f\"{u}r Analysis und Scientific Computing, Technische Universit\"{a}t Wien, Wiedner Hauptstrasse 8-10, A--1040 Wien, Austria.}
 \and 
David W\"{o}rg\"{o}tter\thanks{(david.woergoetter@tuwien.ac.at), Institut f\"{u}r Analysis und Scientific Computing, Technische Universit\"{a}t Wien, Wiedner Hauptstrasse 8-10, A--1040 Wien, Austria.}}

\date{\today}

\begin{document}

\maketitle
\begin{abstract}
		We consider the time-harmonic Maxwell equations with impedance boundary conditions on a bounded Lipschitz domain $\Omega$ with analytic boundary $\Gamma$. We suppose that $\Omega$ consists of multiple subdomains, and that the permeability and permittivity tensors are analytic \revtwo{on} every subdomain, but may jump across subdomain interfaces. Under these conditions we show that for any wavenumber $k\in\Co$ with $|k|\geq 1$ for which Maxwell's equations are polynomially well-posed, a Galerkin discretization based on N\'{e}d\'{e}lec elements of order $p$ on a mesh with mesh width $h$ is quasi-optimal, provided that there holds the wavenumber-explicit scale resolution condition a) that $|k|h/p$ is sufficiently small and b) that $p/\log |k|$ is bounded from below.
\end{abstract}

\section{Introduction}

When solving a \revtwo{wave propagation problem} at a possibly large wavenumber $k>0$ using a \revtwo{low} order Galerkin scheme based on a mesh with mesh width $h$, one usually has to deal with the pollution effect: Even if the number of degrees of freedom per wavelength is held constant for any given wavenumber $k$, the error between the exact solution $\solu$ and the discrete solution $\solu_h$ grows as $k$ increases, \revision{that is,} the gap between the Galerkin error and the best approximation in the discrete space widens as the wavenumber $k$ grows. 

\medskip

In practice, this effect means that under the assumption of a growing wavenumber $k$, asymptotic convergence of a low order Galerkin scheme can only be expected \revtwo{under restrictive conditions on the mesh width $h$ and the wavenumber $k$}. For example, \revtwo{for the Helmholtz equation it can be shown that the corresponding lowest-order finite element method reaches asymptotic convergence only under the condition that $k^2h$ is sufficiently small \cite{BookIhlenburg, HelmholtzIhlenburgBabuska_h, HelmholtzSauter}.} For practical computations, the requirement of $k^2h$ being small is a serious problem since for large $k$, extremely fine meshes are needed to \revtwo{control} pollution. 

During the last three decades numerical methods for the Helmholtz equation such as finite element methods \cite{HelmholtzIhlenburgBabuska_p,HelmholtzDTNMelenk, HelmholtzMelenkSauter, HelmholtzAinsworth,HelmholtzFullspaceSpence, HelmholtzPMLSpence,HelmholtzSpenceGalkowski,HelmholtzBernkopf} or discontinuous Galerkin methods \cite{PaperFengWuDGHelmholtz, PaperFengWuhpDGHelmholtz} have been studied extensively. For finite element methods it turns out that higher order methods are much better suited to solve the Helmholtz equation in the high wavenumber regime than low order methods. In particular, the pollution effect occurring in the numerical simulation of the Helmholtz equation can be greatly reduced by using $hp$-FEM with a sufficiently large local polynomial degree $p$.


\medskip

The natural follow-up question is, whether these results for the Helmholtz equation can be transferred to the time-harmonic Maxwell equations. Compared to the Helmholtz equation, the theory of numerical pollution in Maxwell's equations is far less developed. At least for real-valued and constant coefficients, the pollution effect in $hp$-FEM for Maxwell's equations has been studied in \cite{MaxwellAinsworth} as well as in \cite{MaxwellTransparentMelenk} for transparent boundary conditions and \cite{MaxwellImpedanceMelenk} for impedance boundary conditions. The key result from \cite{MaxwellTransparentMelenk, MaxwellImpedanceMelenk} is that a finite element approximation of \revision{fixed} local polynomial order $p$ based on a mesh with mesh width $h$ of Maxwell's equations at wavenumber $k$ does not suffer from pollution provided that there holds the scale-resolution condition that a) $hp/k$ is sufficiently small and b) $p/\log k$ is bounded from below. 

\medskip 

\revision{The analysis performed in} \cite{MaxwellTransparentMelenk, MaxwellImpedanceMelenk} \revision{restricts to} Maxwell's equations with constant and real-valued coefficients. For non-constant and matrix valued piecewise smooth coefficients, wavenumber-explicit $hp$-finite element convergence results seem to be more sparse. In \cite{MaxwellChaumont} it is shown that for piecewise smooth coefficients, a finite element method based on \ned\ elements is quasi-optimal as long as $hk\leq C_p$ for a constant $C_p>0$ that depends on the geometry, on the coefficients and also on the local polynomial degree $p$ of the finite element scheme.

In the recent work \cite{MaxwellSpence} it is proved that for fixed-order methods the pollution effect is mitigated under the scale-resolution condition that \revtwo{$\rho(k)(kh)^{2p}$ is sufficiently small, where $\rho(k)>0$ is the norm of the Maxwell solution operator cf. \eqref{rhodef} below}. When comparing this condition to the scale-resolution condition in the scalar case \cite{MaxwellTransparentMelenk, MaxwellImpedanceMelenk}, namely a) that $kh/p$ is sufficiently small and b) that $p/\log k$ is bounded from below, it strikes that $\rho(k)(kh)^{2p}$ being sufficiently small is a much more restrictive condition.
The reason for this discrepancy is that the more relaxed scale-resolution condition from \cite{MaxwellTransparentMelenk, MaxwellImpedanceMelenk} crucially relies on the fact that the solution $\solu$ of Maxwell's equations is analytic, provided that the geometry and the given input data are analytic \cite{MaxwellMyself2}. 

\medskip

The aim of this work is to generalize \cite{MaxwellImpedanceMelenk} to the case of piecewise analytic matrix-valued coefficients. That is, we prove that for piecewise analytic and matrix-valued coefficients, the less restrictive scale-resolution condition a) that $kh/p$ is sufficiently small and b) that $p/\log k$ is bounded from below is sufficient to mitigate pollution in $hp$-FEM. In this way, we provide a companion result to \cite{MaxwellSpence}; where the analysis of fixed-order methods in \cite{MaxwellSpence} relies only on piecewise \revtwo{smooth} coefficients and requires the more restrictive scale-resolution condition, we assume that the coefficients are piecewise analytic and show that the less restrictive scale-resolution condition suffices to mitigate pollution in $hp$-FEM.

\bigskip 

For the rest of this work we consider the time-harmonic Maxwell equations with impedance boundary conditions at a wavenumber $k\in\Co$ with\footnote{The condition $|k|\geq 1$ can be replaced by $|k|\geq c_0>0$ at the expense of all involved constants depending also on $c_0$.} $|k|\geq 1$ posed on a bounded \revision{and simply connected} Lipschitz domain $\Omega$ with analytic boundary $\Gamma:=\partial\Omega$. That is, given a right-hand side $\solf$ and impedance boundary data $\solgi$ we look for a solution $\solu$ of 
\begin{equation}\label{Maxwellorig}
\begin{alignedat}{2}
		\curl\mu^{-1}\curl\solu-k^2\varepsilon\solu &= \solf\quad &&{\rm in}\quad \Omega, \\
		\left(\mu^{-1}\curl\solu\right)\times\soln-ik\zeta\solu_T &= \solgi \quad &&{\rm on}\quad  \Gamma,
\end{alignedat}
\end{equation}
where $i:=\sqrt{-1}$ is the imaginary unit, $\soln$ is the outer unit normal to $\Gamma$ and $\solu_T:=\soln\times(\solu\times\soln)$ denotes the tangential trace of $\solu$. Concerning the coefficients we suppose that $\mu^{-1}$ and $\varepsilon$ 
are complex-valued tensor fields which are uniformly coercive and piecewise analytic in $\Omega$, but may be discontinuous across certain mutually disjoint analytic surfaces $\interf_1,\ldots,\interf_r$ contained in $\Omega$. Furthermore, $\zeta:\Gamma\rightarrow\Co^{3\times 3}$ is supposed to be uniformly coercive and analytic on $\Gamma$ and should satisfy $(\zeta\solz)_T = \zeta\solz_T$ for all vector fields $\solz$ on $\Gamma$. \revision{Let us highlight that the ideas and techniques of this work can also be applied in the case of other boundary conditions, e.g., natural boundary conditions that read as $\mu^{-1}\curl\solu\cdot\soln = \solgn$ on $\Gamma$ for some given data $\solgn$ or essential boundary conditions given by $\solu_T=0$ on $\Gamma$. However, in order to not overload this work we consider only impedance conditions; at least for the regularity splitting result provided by Theorem~\ref{Mainresult1} a discussion of other boundary conditions can be found in \cite{MyThesis}.}

\medskip

We notice that our choice of wavenumbers $k\in\Coone:=\{z\in\Co\ |\ |z|\geq 1\}$ is rather untypical; usually one considers only wavenumbers $k>0$. We highlight that Maxwell's equations \eqref{Maxwellorig} are not uniquely solvable for all $k\in\Coone$, this follows from the Fredholm alternative, cf. \cite[Ch.~4]{BookMonk}. In practice, this is not a problem since one is rarely interested in solving \eqref{Maxwellorig} for any $k\in\Coone$. Instead, one is often interested in \eqref{Maxwellorig} only at certain wavenumbers $k\in\setu$, where $\setu\subseteq \Coone$ is a set of \textit{interesting} wavenumbers that is fixed a priori. For example, choices of $\setu$ could be $\setu = \R\cap\Coone$ or $\setu=\{z\in\Coone\ |\ 0<\imagpart z < a\}$ for some constant $a>0$. 

In fact, both of our main results, Theorem~\ref{Mainresult1} and Theorem~\ref{Mainresult2} rely on the assumption that $k\in\setu$ where $\setu\subseteq\Coone$ is such that Maxwell's equations \eqref{Maxwellorig} are polynomially well-posed for all $k\in\setu$, see Assumption~\ref{assumption1} and Assumption~\ref{assumption3} below. Of course all involved constants then also depend on the particular choice of $\setu$, but this is not a problem since $\setu$ is fixed a priori and subsequently varying the specific wavenumber $k\in\setu$ then does not influence the involved constants.

\medskip

The outline of this work is as follows: In Section~\ref{notation} we explain the notation used throughout this work, discuss our assumptions on the underlying geometry and coefficients in more detail and culminate with the statement of the two main results of this work, namely, Theorem~\ref{Mainresult1} and Theorem~\ref{Mainresult2}. The first main result, Theorem~\ref{Mainresult1}, establishes the existence of a regularity splitting of solutions $\solu$ of Maxwell's equations~\eqref{Maxwellorig} akin to \cite[Thm.~7.3]{MaxwellImpedanceMelenk}. The second main result, Theorem~\ref{Mainresult2}, shows that a Galerkin discretization based on \ned-type I elements of order $p$ on a mesh with mesh width $h$ is quasi-optimal provided that there holds the aforementioned scale-resolution condition a) that $|k|h/p$ is sufficiently small and b) that $p/\log |k|$ is bounded from below.

\smallskip

After the statement of our main results we proceed with two auxiliary sections, namely, Section~\ref{freqsec} and Section~\ref{auxmax} in which we provide the tools that are necessary to prove Theorem~\ref{Mainresult1}. To be specific, in Section~\ref{freqsec} we collect and extend some results from \cite[Sec.~6]{MaxwellImpedanceMelenk} concerning frequency splitting operators on vector fields, and in Section~\ref{auxmax} we study a modified version of \eqref{Maxwellorig} where the coefficients $\varepsilon$ and $\zeta$ are chosen such that existence and uniqueness of solutions can be proved via the Lax Milgram lemma. 

\smallskip

The results from Section~\ref{freqsec} and Section~\ref{auxmax} are subsequently combined in Section~\ref{Mainres1proof}, and following the ideas from \cite[Sec.~7]{MaxwellImpedanceMelenk} this leads to a proof of our first main result, Theorem~\ref{Mainresult1}.

\smallskip

After having dealt with Theorem~\ref{Mainresult1}, we turn our attention to our second main result, namely, Theorem~\ref{Mainresult2}. In Section~\ref{discretizationsec} we discuss the Galerkin discretization based on \ned-type I elements of order $p\geq 0$ corresponding to a triangulation with mesh width $h$ and recall certain $hp$-approximation operators from \cite{MaxwellTransparentMelenk, MaxwellImpedanceMelenk} as well as the recently constructed interpolating projection operators from \cite{OperatorsRojik}.

\smallskip

Finally, in Section~\ref{femanalysis} we conclude the proof of our second main result, Theorem~\ref{Mainresult2}. Starting from a recap of the Schatz argument \cite{Schatz} we extend the arguments from \cite[Sec.~9]{MaxwellImpedanceMelenk} and exploit that the error between the exact solution $\solu$ and the finite-element approximation $\solu_h$ is related to the regularity of a dual problem. Using the regularity splitting result from Theorem~\ref{Mainresult1} we can then thoroughly discuss the regularity of this dual problem, which subsequently leads to an elegant proof of Theorem~\ref{Mainresult2}. 

At last, in Section~\ref{numerics} we conclude this work by some numerical experiments that illustrate our findings.

\section{General notation and main results}\label{notation}
For any two vectors $\solw,\solz\in\Co^3$ with $\solw = (w_1, w_2, w_3)^T$ and $\solz=(z_1,z_2,z_3)^T$ we set $\solw\cdot\solz:=\sum_{i=1}^3w_iz_i$ and write $\SCP{\solw}{\solz}{}:=\solw\cdot\overline{\solz}$ for the scalar product between $\solw$ and $\solz$, where $\overline{\solz}:=(\overline{z_1},\overline{z_2},\overline{z_3})^T$ denotes the complex conjugate of $\solz$. 
Furthermore, the cross product between the vectors $\solw$ and $\solz$ is defined in the usual way as $\solw\times\solz:=(w_2z_3-w_3z_2, w_3z_1-w_1z_3, w_1z_2-w_2z_1)^T$, and from the introduction we recall $\Coone:=\{z\in\Co\ |\ |z|\geq 1\}$. For any sufficiently smooth vector field $\solv = (v_1, v_2, v_3)^T$ we define its curl and divergence by
$$
		\curl\solv := \left(\partial_y v_3-\partial_z u_2, \partial_z u_1-\partial_x u_3, \partial_x u_2-\partial_y u_1\right)^T\quad\text{and}\quad \diverg\solv:=\partial_x v_1+\partial_y v_2+\partial_z v_3. $$

 As usual, for $\Omega\subseteq\R^3$ let $\SL(\Omega)$ denote the Lebesgue space of complex-valued square integrable functions, and define its vector-valued version as $\VSL(\Omega):=(\SL(\Omega))^3$. Furthermore, $\SLinf(\Omega)$ is the space of all essentially bounded functions on $\Omega$ and $\Linfty(\Omega):=\left(\SLinf(\Omega)\right)^3$.
  For a bounded Lipschitz domain $\Omega\subseteq \R^3$ and $\ell\in\N_0$, the space $\SHM{\ell}(\Omega)$ is the usual Sobolev space of order $\ell$, see \cite[Ch.~3]{BookMcLean}, and $\SHM{\ell}_0(\Omega)$ denotes the closure of $\CM{\infty}_0(\Omega)$ in $\SHM{\ell}(\Omega)$. In order to deal with vector fields, we define the vector-valued space $\VSHM{\ell}(\Omega):=(\SHM{\ell}(\Omega))^3$.
  In addition, we define the spaces
  \begin{align*}
  \Hcurl := \{\solu\in\VSL(\Omega)\ |\ \curl \solu\in\VSL(\Omega)\}\quad {\rm and}\quad \Hdiv:=\{\solu\in\VSL(\Omega)\ |\ \diverg \solu\in\SL(\Omega)\}.
  \end{align*}

\medskip

Furthermore, for a sufficiently smooth, \revision{closed and orientable} surface $\Sigma$ and $s> 0$, let $\SHM{s}(\Sigma)$ be the fractional Sobolev space of index $s$ with dual space $\SHM{-s}(\Sigma)$, see \cite[Ch.~3]{BookMcLean}, and let $\VSHM{s}(\Sigma):=(\SHM{s}(\Sigma))^3$. Moreover, the space of square-integrable tangent fields is given by
\begin{align*}
\VSL_T(\Sigma):=\{\solv\in\VSL(\Sigma)\ |\ \solv\cdot\soln = 0\},
\end{align*}
where $\soln$ is the outer unit normal to $\Sigma$, and for $s\geq 0$ we set 
\begin{align*}
\VSHM{s}_T(\Sigma):=\VSL_T(\Sigma)\cap\VSHM{s}(\Sigma),\ \ {\rm as\ well\ as}\ \ \VSHM{-s}_T(\Sigma):=(\VSHM{s}_T(\Sigma))'.
\end{align*}

For a bounded Lipschitz domain $\Omega\subseteq\R^3$ we say that a function $u$ or a vector field $\solu$ is smooth in $\Omega$ if $u\in\CM{\infty}(\Omega)$ or $\solu\in\left(\CM{\infty}(\Omega)\right)^3$, respectively. In addition, a function $u$ or a vector field $\solu$ is said to be smooth up to the boundary of $\Omega$ if $u$ or $\solu$ can be extended to a function $\widetilde{u}\in\CM{\infty}(\R^3)$ or a vector field $\widetilde{\solu}\in\left(\CM{\infty}(\R^3)\right)^3$, respectively. In this case we also write $u\in\CM{\infty}(\overline{\Omega})$ for functions $u$ and $\solu\in\left(\CM{\infty}(\overline{\Omega})\right)^3$ for vector fields $\solu$.

\medskip

Henceforth, we say that a tensor field $\nu:\Omega\rightarrow\Co^{3\times 3}$ or $\nu:\Gamma\rightarrow\Co^{3\times 3}$ is uniformly coercive if there exists a constant $c>0$ and a complex number $\alpha$ with $|\alpha|=1$ such that
\begin{align}\label{coercivitydef}
		\revision{\forall\solz\in\Co^3:\ }	\realpart\SCP{\alpha\nu\solz}{\solz}{}\geq c\norm{\solz}{}^2
\end{align}
uniformly in $\Omega$ or uniformly on $\Gamma$, respectively.

\medskip

The following definition clarifies the notion of functions and vector fields that are analytic up to the boundary of a bounded Lipschitz domain $\Omega$. For a multiindex $\alpha\in\N_0^3$ and a vector field $\solv = (v_1, v_2, v_3)^T$ we set $\D^{\alpha}(\solv):=\left(\D^{\alpha}(v_1),\D^{\alpha}(v_2), \D^{\alpha}(v_3)\right)^T$. 

\begin{definition}\label{analytictensors}
		For a bounded Lipschitz domain $\Omega\subseteq\R^3$, a wavenumber $k\in\Coone$ and given constants $\omega\geq 0,M>0$, the sets $\anaf{\Omega}$ and $\anavec{\Omega}$ consist of all functions $u:\Omega\rightarrow\Co$ and vector fields $\solu:\Omega\rightarrow\Co^3$ which satisfy 
\begin{align*}
		\revision{\forall\ell\in\N_0:\ }	\sum_{|\alpha|=\ell}\norm{\D^{\alpha}(u)}{\SL(\Omega)}\leq \omega M^{\ell}(\ell+|k|)^{\ell} \quad {\rm and}\quad 
		\revision{\forall\ell\in\N_0:\ }\sum_{|\alpha|=\ell}\norm{\D^{\alpha}(\solu)}{\VSL(\Omega)}\leq \omega M^{\ell}(\ell+|k|)^{\ell}.
\end{align*}

\smallskip

Moreover, the symbol $\anatens{\Omega}$ denotes the space of all tensor fields $\nu:\Omega\rightarrow\Co^{3\times 3}$ for which there exist constants $\omega\geq 0, M>0$ such that 
\begin{align*}
		\revision{\forall\ell\in\N_0:\ }\sum_{i,j=1}^3\sum_{|\alpha|=\ell}\norm{\D^{\alpha}(\nu_{i,j})}{\SL(\Omega)}\leq \omega M^{\ell}\ell^{\ell}.
\end{align*}
\end{definition}

\begin{remark}
		The reader will immediately notice that the definitions of $\anaf{\Omega}$ and $\anavec{\Omega}$ explicitly include $\omega, M$ and $k$, whereas the definition of $\anatens{\Omega}$ does not. This is due to the fact that the proof of Theorem~\ref{Mainresult1} requires a meticulous tracking of the growth of the derivatives of wavenumber-dependent vector fields. The coefficients $\mu^{-1}$, $\varepsilon$ and $\zeta$, however, are not wavenumber-dependent, and the growth of their derivatives is of less importance for this work, hence we suppress $\omega$ and $M$ and leave away the wavenumber-dependency in the notion of the analyticity classes $\anatens{\Omega}$.
\end{remark}

While Definitions~\ref{analytictensors} clarifies the notion of analytic functions, vector- and tensor fields on $\Omega$, the following definition introduces analyticity classes on the boundary of a bounded Lipschitz domain with analytic boundary. 

\begin{definition}\label{assumptionimpedance}
		Let $\Omega\subseteq\R^3$ be a bounded Lipschitz domain with analytic boundary $\Gamma$ and suppose that constants $\omega\geq 0,M>0$ as well as $k\in\Coone$ are given. The sets $\anagammaf$ and $\anagammavec{\Gamma}$ consist of all functions $g:\Gamma\rightarrow\Co$ and tangent fields $\solg_T:\Gamma\rightarrow\Co^3$ for which there exist continuations $g^*$ and $\solg^*$ to a tubular neighbourhood $\mathcal{U}$ of $\Gamma$ which satisfy
		\begin{align*}
				\revision{\forall\ell\in\N_0:\ }\sum_{|\alpha|=\ell}\norm{\D^{\alpha}(g^*)}{\SL(\mathcal{U})}\leq \omega M^{\ell}(\ell+|k|)^{\ell}\quad {\rm and}\quad 
				\revision{\forall\ell\in\N_0:\ }\sum_{|\alpha|=\ell}\norm{\D^{\alpha}(\solg^*)}{\VSL(\mathcal{U})}\leq \omega M^{\ell}(\ell+|k|)^{\ell}.
		\end{align*}

		Similarly, the space $\anagamma$ consists of all tensor fields $\nu:\Gamma\rightarrow\Co^{3\times 3}$ for \revision{which} there exist constants $\omega\geq 0, M>0$ and an extension $\nu^*$ to a tubular neighbourhood $\mathcal{U}$ of $\Gamma$ such that
\begin{align*}
		\revision{\forall\ell\in\N_0:\ }\sum_{i,j=1}^3\sum_{|\alpha|=\ell}\norm{\D^{\alpha}(\nu^*_{i,j})}{\SL(\mathcal{U})}\leq \omega M^{\ell}\ell^{\ell}.	
\end{align*}
\end{definition}

%
%
%
%

\subsection{$\ana$-partitions and classes of piecewise analytic functions}

Let $\Omega\subseteq\R^3$ be a bounded Lipschitz domain with boundary $\Gamma$. 
Throughout this work we suppose that $\Omega$ is partitioned into subdomains $\Gp_1,\ldots,\Gp_n$ such that the boundaries of these subdomains form a family of mutually disjoint closed\footnote{By {\it closed surface} we denote a compact surface without boundary.} and analytic surfaces inside of $\Omega$. 
The following definition makes this precise.

\begin{definition}\label{partitiondef}
		An $\ana$-partition $\Gp$ is a tuple $\Gp = \geom$ which consists of domains $\Omega, \Gp_1,\ldots, \Gp_n\subseteq\R^3$ satisfying 
		\begin{itemize}
				\item[(i)] The domains $\Omega, \Gp_1,\ldots,\Gp_n$ are bounded Lipschitz domains with $\Omega$ being simply connected and with the domains $\Gp_1,\ldots,\Gp_n$ being mutually disjoint and satisfying $\overline{\Omega} = \overline{\Gp_1}\cup\ldots\cup\overline{\Gp_n}$.
				\item[(ii)] The boundary $\Gamma:=\partial\Omega$ is simply connected and analytic, and the subdomain boundaries $\partial\Gp_1,\ldots,\partial\Gp_n$ are analytic. 
				\item[(iii)] If $n>1$, there exist bounded Lipschitz domains $V_1,\ldots,V_r$ each having a simply connected analytic boundary $\interf_j:=\partial V_j$ such that $\Gamma, \interf_1,\ldots,\interf_r$ are mutually disjoint and
\begin{align*}
		\Gamma\cup\bigcup_{j=1}^r\interf_j = \bigcup_{j=1}^n\partial\Gp_j.
\end{align*}
The surfaces $\interf_1,\ldots,\interf_r$ are called subdomain interfaces and their union $\interf := \interf_1\cup\ldots\cup\interf_r$ is referred to as the union of subdomain interfaces.
\end{itemize}
\end{definition}

The point of introducing subdomains in the preceding fashion is to incorporate the location of (possible) discontinuities of piecewise regular coefficients $\mu^{-1}$ and $\varepsilon$ into the geometry of the problem. If $n=1$, then $\Omega$ has only one subdomain, namely itself, hence in this case there are no subdomain interfaces present.

\begin{remark}
We highlight that requirement (iii) in Definition~\ref{partitiondef} implies that every connected component of $\partial\Gp_i$ either coincides with $\Gamma$ or with a subdomain interface $\interf_j$. Moreover, for every $\interf_j$ there are precisely two subdomains $\Gp_i$ and $\Gp_h$ such that $\interf_j = \partial\Gp_i\cap\partial\Gp_h$. 
In particular, there may be no point in $\Omega$ where three or more subdomains meet.

Moreover, the assumptions that $\Omega$, $\Gamma$ and all interface components are simply connected simplify some of our proofs but they are not essential for our first main result, Theorem~\ref{Mainresult1}. In fact, Theorem~\ref{Mainresult1} remains valid if $\Omega, \Gamma$ or any of the subdomain interfaces $\interf_j$ are not simply connected with only minor modifications necessary cf. \cite[Ch.9]{MyThesis}. Our second main result, Theorem~\ref{Mainresult2}, however, needs a simply connected domain $\Omega$ because its proof crucially relies on a (discrete) sequence property of the de Rham complex over $\Omega$.
\end{remark}

The notion of $\ana$-partitions is the basis for broken Sobolev spaces, which play an important role in this work. 
Let $\Gp=\geom$ be an $\ana$-partition and assume $m\in\N_0$. We define
\begin{align*}
\PSHM{m}(\Gp):=\{u\in\SL(\Omega)\ \big|\ u\vert_{\Gp_i}\in\SHM{m}(\Gp_i)\ {\rm for}\ i=1,\ldots,n\} \quad{\rm and}\quad \PVSHM{m}(\Gp):=(\PSHM{m}(\Gp))^3,
\end{align*}
and for $u\in\PSHM{m}(\Gp)$ and $\solu=(u_1,u_2,u_3)^T\in\PVSHM{m}(\Gp)$ we set
\begin{align*}
		\norm{u}{\PSHM{m}(\Gp)}^2:=\sum_{j=1}^n\norm{u}{\SHM{m}(\Gp_j)}^2\quad{\rm and} \quad \norm{\solu}{\PVSHM{m}(\Gp)}^2:=\sum_{i=1}^3\norm{u_i}{\PSHM{m}(\Gp)}^2.
\end{align*}
Furthermore, for $m\in\N_0$ we introduce the spaces of vector fields with piecewise regular curl and divergence
\begin{align*}
\PVHcurl{m}:=\{\solv\in\Hcurl\ \big|\ \curl\solv\vert_{\Gp_i}\in\VSHM{m}(\Gp_i)\ {\rm for}\ i=1,\ldots,n\}
\end{align*}
and
\begin{align*}
\PVHdiv{m}:=\{\solv\in\Hdiv\ \big|\ \diverg\solv\vert_{\Gp_i}\in\SHM{m}(\Gp_i)\ {\rm for}\ i=1,\ldots,n\}.
\end{align*}

For the rest of this work we call functions or vector fields piecewise smooth on an $\ana$-partition $\Gp$ if they are smooth up to the boundary of every subdomain. The following definition makes this precise.
\begin{definition}\label{pwsmooth}
		Let $\Gp=\geom$ be an $\ana$-partition. A function $v:\Omega\rightarrow\Co$ or a vector field $\solv:\Omega\rightarrow\Co^3$ is called piecewise smooth if for all $\ell\in\N_0$ there holds $v\in\PSHM{\ell}(\Gp)$ or $\solv\in\PVSHM{\ell}(\Gp)$, respectively.
\end{definition}
We notice that the terminology {\it piecewise smooth} is justified by the Sobolev embedding theorem, according to which a piecewise smooth function or vector field is indeed smooth up to the boundary of every subdomain.

\medskip

We adapt Definition~\ref{analytictensors} to the setting of $\ana$-partitions.

\begin{definition}\label{pwanalyticdef}
		Let $\Gp=\geom$ be an $\ana$-partition and suppose that a wavenumber $k\in\Coone$ and constants $\omega\geq 0, M>0$ are given. The sets $\anapw{\Gp}$ and $\anavecpw{\Gp}$ consist of all functions $u:\Omega\rightarrow\Co$ and vector fields $\solu:\Omega\rightarrow\Co^3$ which satisfy 
\begin{align*}
		\forall\ell\in\N_0:\ \sum_{i=1}^n\sum_{|\alpha|=\ell}\norm{\D^{\alpha}(u)}{\SL(\Gp_i)}\leq \omega M^{\ell} (\ell+|k|)^{\ell}\quad {\rm and}\quad 
		\forall\ell\in\N_0:\ \sum_{i=1}^n\sum_{|\alpha|=\ell}\norm{\D^{\alpha}(\solu)}{\VSL(\Gp_i)}\leq \omega M^{\ell} (\ell+|k|)^{\ell}.
\end{align*}

Similarly, the space $\anatenspw{\Gp}$ consists of all tensor fields $\nu:\Omega\rightarrow\Co^{3\times 3}$ for which there exist constants $\omega\geq 0, M>0$ such that
\begin{align*}
		\forall\ell\in\N_0:\ \sum_{i=1}^n\sum_{j,m=1}^3\norm{\D^{\alpha}(\nu_{j,m})}{\SL(\Gp_i)}\leq \omega M^{\ell}\ell^{\ell}.
\end{align*}
\end{definition}

For the rest of this work we will always work with coefficients $\mu^{-1}, \varepsilon\in\anatenspw{\Gp}$ and $\zeta\in\anatens{\Gamma}$. In addition to that, we will always assume that the coefficients $\mu^{-1},\varepsilon$ are uniformly coercive on $\Omega$ and that $\zeta$ is uniformly coercive on $\Gamma$ and well-behaved when applied to tangent fields, see Assumption~\ref{assumption1} below.

\subsection{Traces of $\Hcurl$ and variational formulation of Maxwell's equations}

Let us assume that $\Omega$ is a bounded Lipschitz domain with smooth boundary $\Gamma$. With $\curl_{\Gamma}$ and $\diverg_{\Gamma}$ denoting the surface curl and surface divergence from e.g. \cite{BookMonk, BookNedelec, MaxwellMyself, MaxwellImpedanceMelenk} we define the auxiliary spaces 
\begin{align*}
		\Hcurlgamma{\Gamma}&:=\left\{\solu\in\VSHM{-1/2}_T(\Gamma)\ \big|\ \curl_{\Gamma}\solu\in\SHM{-1/2}(\Gamma)\right\}, \\
		\Hdivgamma{\Gamma}&:=\left\{\solu\in\VSHM{-1/2}_T(\Gamma)\ \big|\ \diverg_{\Gamma}\solu\in\SHM{-1/2}(\Gamma)\right\},
\end{align*}
and equip them with norms
\begin{align*}
		\norm{\solu}{\Hcurlgamma{\Gamma}}^2&:=\norm{\solu}{\VSHM{-1/2}_T(\Gamma)}^2+\norm{\curl_{\Gamma}\solu}{\SHM{-1/2}(\Gamma)}^2,\\
		\norm{\solu}{\Hdivgamma{\Gamma}}^2&:=\norm{\solu}{\VSHM{-1/2}_T(\Gamma)}^2+\norm{\diverg_{\Gamma}\solu}{\SHM{-1/2}(\Gamma)}^2.
\end{align*}
There holds the following trace result, see \cite[Thm.~3.29]{BookMonk} or \cite[Thm.~5.4.2]{BookNedelec}.

\begin{proposition}\label{traceprop}
		Let $\Omega\subseteq\R^3$ be a bounded Lipschitz domain with boundary $\Gamma$ and outer unit normal $\soln$. 
		We consider the maps $\Pi_{T}$ and $\Pi_{t}$, which for vector fields $\solv\in\left(\CM{\infty}(\overline{\Omega})\right)^3$ are defined as 
	\begin{align*}
			\Pi_{T}\solv := \soln\times(\solv\vert_{\Gamma}\times\soln)\quad {\rm and}\quad \Pi_{t}\solv := \solv\vert_{\Gamma}\times\soln.
	\end{align*}
	These maps extend to bounded  and surjective operators $\Pi_{T}:\Hcurl\rightarrow\Hcurlgamma{\Gamma}$ and $\Pi_{t}:\Hcurl\rightarrow\Hdivgamma{\Gamma}$. Moreover, there exist bounded lifting operators $\liftcurl:\Hcurlgamma{\Gamma}\rightarrow\Hcurl$ and $\liftdiv:\Hdivgamma{\Gamma}\rightarrow\Hcurl$. 
\end{proposition}

\begin{remark}
In order to shorten notation we we will often abbreviate $\Pi_T\solv$ by $\solv_T$ and $\Pi_t\solv$ by $\solv_t$, respectively. 
\end{remark}

At this point we have all necessary tools to clarify the notion of a {\it weak solution} of Maxwell's equations \eqref{Maxwellorig} at a wavenumber $k\in\Coone$. We assume that coefficients $\mu^{-1}, \varepsilon\in\anatenspw{\Gp}$ and $\zeta\in\anatens{\Gamma}$ are given. To formulate the variational formulation of \eqref{Maxwellorig} we define the energy space
\begin{align}\label{energyspace}
			\HXI:=\{\solu\in\Hcurl\ |\ \solu_T\in\VSL(\Gamma)\} 
\end{align}
together with its companion space of scalar potentials
\begin{align*}
		\Honeimp=\{\varphi\in\SHM{1}(\Omega)\ |\ \varphi\vert_{\Gamma}\in\SHM{1}(\Gamma)\}.
\end{align*}
We equip $\HXI$ with the wavenumber-dependent norm
\begin{align*}
		\norm{\solu}{\HXIK}^2:=\norm{\curl\solu}{\VSL(\Omega)}^2+|k|^2\norm{\solu}{\VSL(\Omega)}^2+|k|\norm{\solu_T}{\VSL_T(\Gamma)}^2.
\end{align*}

With these definitions, let a right-hand side $\solf\in\VSL(\Omega)$ and impedance data $\solgi\in\VSL_T(\Gamma)$ be given. Then, a vector field $\solu\in\HXI$ is called a weak solution of Maxwell's equations \eqref{Maxwellorig} if
\begin{align}\label{weakformulation}
		\forall\solv\in\HXI:\ \bfa(\solu,\solv)
					= \SCP{\solf}{\solv}{\VSL(\Omega)}+\SCP{\solgi}{\solv_T}{\VSL_T(\Gamma)},
\end{align}
where
\begin{align}\label{Adef}
\bfa(\solu,\solv):= \SCP{\mu^{-1}\curl \solu}{\curl\solv}{\VSL(\Omega)}-k^2\SCP{\varepsilon\solu}{\solv}{\VSL(\Omega)}-ik\SCP{\zeta\solu_T}{\solv_T}{\VSL_T(\Gamma)}
\end{align}
for all $\solu,\solv\in\HXI$.

\bigskip 

Establishing conditions on the coefficients $\mu^{-1}, \varepsilon, \zeta$ and the wavenumber $k$ which guarantee unique solvability of \eqref{weakformulation} is beyond the scope of this work. Even in the homogeneous case $\mu^{-1}=\varepsilon=\zeta=1$, there always exists a discrete set of wavenumbers such that \eqref{weakformulation} is not uniquely solvable; this can be proved by showing that the underlying Maxwell operator is Fredholm \cite[Sec.~4.5]{BookMonk}. 

The fact that \eqref{Maxwellorig} is not uniquely solvable for all choices of $k$ is not so much a hindrance in practice as one usually does not attempt to solve Maxwell's equations for all $k\in\Coone$. Instead, one is usually interested in solutions at wavenumbers $k\in\setu$ for some region $\setu\subseteq\Coone$, for example $\setu = \R$ or $\setu$ being the set of all complex numbers with positive but ``small'' imaginary part. 
Therefore it is reasonable to restrict ourselves to the case $k\in\setu\subseteq\Coone$, where $\setu$ is such that unique solvability of \eqref{weakformulation} is guaranteed, and where the given coefficients exhibit sufficiently nice properties. \revtwo{Assumption~\ref{assumption1} below makes these requirements precise.}

\begin{assumption}\label{assumption1}
		Let $\Gp=\geom$ be a given $\ana$-partition and let $\mu^{-1}, \varepsilon\in\anatenspw{\Gp}$ as well as $\zeta\in\anatens{\Gamma}$ be given coefficients.

		We assume that there exists an ellipticity constant $\cco>0$ and complex numbers $\alpha_{\varepsilon}$ and $\alpha_{\zeta}$ satisfying $|\alpha_{\varepsilon}|=|\alpha_{\zeta}|=1$ such that 
		\begin{subequations}
				\begin{alignat}{2}
						&\revision{\forall\solw\in\Co^3:\ }\realpart\SCP{\mu^{-1}\solw}{\solw}{}&&\geq \cco\norm{\solw}{}^2, \label{coercivemu} \\
						&\revision{\forall\solw\in\Co^3:\ }\realpart\SCP{\alpha_{\varepsilon}\varepsilon\solw}{\solw}{}&&\geq \cco\norm{\solw}{}^2, \label{coerciveeps} \\
						&\revision{\forall\solw\in\Co^3:\ }\realpart\SCP{\alpha_{\zeta}\varepsilon\solw}{\solw}{}&&\geq \cco\norm{\solw}{}^2 \label{coercivezeta}
				\end{alignat}
		\end{subequations}
		uniformly on $\Omega$ and $\Gamma$, respectively. 
		In addition we suppose that $\zeta$ satisfies
		\begin{align}\label{zetaproperty}
			(\zeta\solz)_T=\zeta\solz_T
		\end{align}
		for all smooth vector fields $\solz:\Gamma\rightarrow\Co^3$.

		Finally, we assume that $\setu\subseteq\Coone$ is such that for all $k\in\setu$, Maxwell's equations \eqref{Maxwellorig} are well-posed at wavenumber $k$. That is, we assume that for every $k\in\setu$, every $\solf\in\VSL(\Omega)$ and every $\solgi\in\VSL_T(\Gamma)$ there exists a unique weak solution $\solu\in\HXI$ of Maxwell's equations \eqref{Maxwellorig}, and
\begin{align}\label{rhodef}
		\norm{\solu}{\HXIK}\leq \rho(k)\left(\norm{\solf}{\VSL(\Omega)}+\norm{\solgi}{\VSL_T(\Gamma)}\right)
\end{align}
for a number $\rho(k)>0$ independent of $\solf$ and $\solgi$.
\end{assumption}

We mention that \eqref{coercivemu}-\eqref{coercivezeta} is equivalent to the coefficients $\mu^{-1}, \varepsilon$ and $\zeta$ in Maxwell's equations \eqref{Maxwellorig} being uniformly coercive in the sense of \eqref{coercivitydef}, and $\eqref{zetaproperty}$ means that $\zeta$ maps tangent fields again to tangent fields and that $\zeta\soln = \oldlambda\soln$ for a scalar-valued function $\oldlambda$ on $\Gamma$. Due to $\zeta\in\anatens{\Gamma}$, the function $\oldlambda$ is necessarily analytic on $\Gamma$. 

\medskip

Furthermore, we notice that the uniqueness assumption in Assumption~\ref{assumption1} and \eqref{rhodef} assert that the solution map $\solimp:\VSL(\Omega)\times\VSL_T(\Gamma)\rightarrow\HXI$ which maps $\solf,\solgi$ to the unique weak solution $\solu\in\HXI$ is well-defined and bounded with norm $\rho(k)$.

\begin{remark}
		Let us mention that Assumption~\ref{assumption1} is certainly satisfied if $\mu^{-1},\varepsilon\in\anatenspw{\Gp}$ are real-valued, symmetric and uniformly positive definite, the impedance $\zeta$ is a real-valued, analytic and positive function on $\Gamma$ and $\setu\subseteq\Co^+:=\{z\in\Co\ |\ \imagpart z>0\}$. To prove this fact, one can for example use that the underlying Maxwell operator is Fredholm \cite[Sec.~4.5]{BookMonk} together with elementary properties of $\bfa(\cdot,\cdot)$. 

		\smallskip

Moreover, if $\mu^{-1},\varepsilon\in\anatenspw{\Gp}$ are real-valued, symmetric and positive definite, and $\zeta$ is a non-zero real constant, Assumption~\ref{assumption1} holds for all $\setu\subseteq\{z\in\Co\ |\ \imagpart z\geq 0\}$. This follows from the Fredholm property of the underlying Maxwell operator and a unique continuation principle \cite[Thm.~2]{Ball}.
\end{remark}

When conducting our finite element analysis in Section~\ref{femanalysis} below it turns out that the coercivity assumptions \eqref{coercivemu}-\eqref{coercivezeta} are too weak for our purposes. In fact we need some sort of \revtwo{uniform} coercivity \revtwo{(in $\Omega$ or $\Gamma$)} of $\varepsilon$ and $\zeta$, hence from Section~\ref{femanalysis} on we rely on the following assumption:

\begin{assumption}\label{assumption2}
		Let $\Gp=\geom$ be a given $\ana$-partition, let $\mu^{-1}, \varepsilon\in\anatenspw{\Gp}$ as well as $\zeta\in\anatens{\Gamma}$ be given coefficients and let $\setu\subseteq\Coone$ be a region of wavenumbers for which we are interested in the solution of Maxwell's equations \eqref{Maxwellorig}. 

		Under these conditions we assume that $\mu^{-1}$ satisfies \eqref{coercivemu} for an ellipticity constant $\cco>0$, and that $\zeta$ satisfies \eqref{zetaproperty}. Furthermore, we suppose that for every $k\in\setu$ there is a complex number $\alpha_k$ with $|\alpha_k|=1$ such that
		\begin{subequations}
				\begin{alignat}{2}
						&\revision{\forall\solw\in\Co^3:\ }\realpart\alpha_kk^2\SCP{\varepsilon\solw}{\solw}{}&&\geq \cco |k|^2\norm{\solw}{}^2, \label{ezcoerceps} \\
						&\revision{\forall\solw\in\Co^3:\ }\realpart\alpha_kik\SCP{\zeta\solw}{\solw}{}&&\geq \cco |k|\norm{\solw}{}^2 \label{ezcoerczeta}
		\end{alignat}
\end{subequations}
		uniformly on $\Omega$ and on $\Gamma$, respectively.
\end{assumption} 

We notice that \eqref{ezcoerceps}-\eqref{ezcoerczeta} are indeed stronger than \eqref{coerciveeps}-\eqref{coercivezeta}. A tedious calculation shows that Assumption~\ref{assumption2} is certainly satisfied for real-valued, symmetric and uniformly positive definite tensor fields $\mu^{-1},\varepsilon\in\anatenspw{\Gp}$ and $\zeta\in\anatens{\Gamma}$, and \revision{$\setu = \left\{k\in\Co\setminus\{0\}\ |\ \imagpart k\geq 0\right\}$. }

\medskip

The big advantage of Assumption~\ref{assumption2} is that \revision{with $\cco>0$ and $\alpha_k$ being the constants from Assumption~\ref{assumption2}} we can estimate
\begin{align}\label{weakcoercivity}
		\begin{split}
			\norm{\solu}{\HXIK}^2&\leq \cco\realpart\left[\bfa(\solu,\solu)+(1+\alpha_k)\left(k^2\SCP{\varepsilon\solu}{\solu}{\VSL(\Omega)}+ik\SCP{\zeta\solu_T}{\solu_T}{\VSL_T(\Gamma)}\right)\right] \\
							 &\leq \cco\left|\bfa(\solu,\solu)+(1+\alpha_k)\left(k^2\SCP{\varepsilon\solu}{\solu}{\VSL(\Omega)}+ik\SCP{\zeta\solu_T}{\solu_T}{\VSL_T(\Gamma)}\right)\right| 
		\end{split}
	\end{align}
	for all $\solu\in\HXI$ and all $k\in\setu$; an inequality which turns out to be crucial in Section~\ref{Schatzsection}. We highlight that for \eqref{weakcoercivity} it is important that the constant $\alpha_k$ is the same in \eqref{ezcoerceps} and \eqref{ezcoerczeta}.

\bigskip 

While Assumption~\ref{assumption2} strengthens the coercivity assumptions on $\mu^{-1},\varepsilon$ and $\zeta$ compared to Assumption~\ref{assumption1}, it does not assume anything about well-posedness of Maxwell's equations~\eqref{Maxwellorig} at wavenumbers $k\in\setu$. If we believe in Assumption~\ref{assumption1} then for every $k\in\setu$ the weak solution $\solu\in\HXI$ of \eqref{Maxwellorig} exists, is unique and satisfies the a priori estimate \eqref{rhodef}. However, the drawback of \eqref{rhodef} is that the dependence of $\rho(k)$ on $k$ is in general unknown. 

At least under certain conditions on the coefficients $\mu^{-1}, \varepsilon, \zeta$ and the underlying geometry, it can be shown that $\rho(k)$ can be bounded algebraically in $k$ for all $k>0$, see e.g. \cite{MaxwellHiptmair, PaperChaumontMoiolaSpence} or \cite[Sec.~3]{MaxwellImpedanceMelenk}. This motivates the following assumption on which we will rely from Section~\ref{femanalysis} onwards.

\medskip

\begin{assumption}\label{assumption3}
		Let $\Gp=\geom$ be a given $\ana$-partition, let $\mu^{-1}, \varepsilon\in\anatenspw{\Gp}$ as well as $\zeta\in\anatens{\Gamma}$ be given coefficients and let $\setu\subseteq\Coone$ be a region of wavenumbers for which we are interested in the solution of Maxwell's equations \eqref{Maxwellorig}. 

		Under these conditions we assume that there are constants $C_{\setu}>0$ and \revision{$\theta\geq 0$} such that for every $k\in\setu$, every $\solf\in\VSL(\Omega)$ and every $\solgi\in\VSL_T(\Gamma)$ there exists a unique weak solution $\solu\in\HXI$ of Maxwell's equations~\eqref{Maxwellorig}, and  
\begin{align}\label{rhoalgebra}
		\norm{\solu}{\HXIK}\leq C_{\setu}|k|^{\theta}\left(\norm{\solf}{\VSL(\Omega)}+\norm{\solgi}{\VSL_T(\Gamma)}\right).
\end{align}
That is, we assume that $\rho(k)$ from \eqref{rhodef} can be bounded by $\rho(k)\leq C_{\setu}|k|^{\theta}$ for all $k\in\setu$.
\end{assumption}

If \eqref{rhoalgebra} holds it means that the Maxwell problem \eqref{Maxwellorig} is polynomially well-posed in the sense that it is uniquely solvable and the solution depends continuously on the data with a continuity constant that grows at most polynomially in $k$.
It is clear that together, Assumption~\ref{assumption2} and Assumption~\ref{assumption3} are stronger than Assumption~\ref{assumption1}. 

\medskip

 For given coefficients $\mu^{-1}, \varepsilon$ and $\zeta$, let $\mu^{-H}, \varepsilon^H$ and $\zeta^H$ denote their hermitian transposes. Furthermore, for $\setu\subseteq\Coone$ we set $-\overline{\setu}:=\{z\in\Co\ |\ -\overline{z}\in\setu\}$. Some straightforward manipulations then show the following:

 \begin{lemma}
		 Let $\Gp=\geom$ be a given $\ana$-partition and suppose that $\mu^{-1}, \varepsilon\in\anatenspw{\Gp}$ and $\zeta\in\anatens{\Gamma}$ as well as $\setu\subseteq\Coone$ are such that Assumption~\ref{assumption1}, Assumption~\ref{assumption2}, or Assumption~\ref{assumption3} holds. Then, Assumption~\ref{assumption1}, Assumption~\ref{assumption2} or Assumption~\ref{assumption3} is also satisfied by $\mu^{-H}, \varepsilon^H, \zeta^H$ and the region $-\overline{\setu}$, respectively.
 \end{lemma}

In particular, the above lemma asserts that if Assumption~\ref{assumption2} holds, then for every $k\in\setu$ there exists a unimodular number $\alpha_k\in\Co$ such that
\begin{align}\label{weakcoercivitydual}
		\norm{\solu}{\HXIK}^2\leq \cco\left|\bfadual(\solu,\solu)+(1+\alpha_k)\left(\overline{k}^2\SCP{\varepsilon^H\solu}{\solu}{\VSL(\Omega)}-i\overline{k}\SCP{\zeta^H\solu_T}{\solu_T}{\VSL_T(\Gamma)}\right)\right|
\end{align}
for all $\solu\in\HXI$, where $\cco>0$ is independent of $k$ and 
\begin{align}\label{Adualdef}
		\bfadual(\solu, \solv) := \SCP{\mu^{-H}\curl \solu}{\curl\solv}{\VSL(\Omega)}-\overline{k}^2\SCP{\varepsilon^H\solu}{\solv}{\VSL(\Omega)}+i\overline{k}\SCP{\zeta^H\solu_T}{\solv_T}{\VSL_T(\Gamma)}
\end{align}
for all $\solu,\solv\in\HXI$.

\subsection{Main results}\label{mainresultssec}

The first main result of this work asserts that under a little bit more regularity on the data $\solf$ and $\solgi$, the weak solution $\solu$ of Maxwell's equations \eqref{Maxwellorig} can be written as $\solu=\soluH+\soluA+k^{-2}\nabla\psi_{\solf}+k^{-1}\nabla\psi_{\solgi}$ for a $\PVSHM{2}(\Gp)$-regular part $\soluH$ that can be bounded uniformly in $k$, a piecewise analytic part $\soluA$ for which $k$-explicit bounds are available and two auxiliary gradient fields $\nabla\psi_{\solf}$ and $\nabla\psi_{\solgi}$. 
Besides being interesting in itself, this decomposition is fundamental for the proof of our second main result, Theorem~\ref{Mainresult2}.

\begin{theorem}\label{Mainresult1}
		Let $\Gp=\geom$ be an $\ana$-partition and let $\mu^{-1},\varepsilon\in\anatenspw{\Gp}$, $\zeta\in\anatens{\Gamma}$ and $\setu\subseteq\Coone$ be such that Assumption~\ref{assumption1} holds. 		
		For any wavenumber $k\in\setu$ let $\solu$ be the weak solution of Maxwell's equations \eqref{Maxwellorig} at wavenumber $k$ corresponding to a right-hand side $\solf\in\Hdiv\cap\PVSHM{1}(\Gp)$ and impedance data $\solgi\in\VSHM{1/2}_T(\Gamma)$. 

Then, $\solu$ can be written as $\solu = \soluH+\soluA+k^{-2}\nabla\psi_{\solf}+k^{-1}\nabla\psi_{\solgi}$, and $\soluH$, $\psi_{\solf}$ and $\psi_{\solgi}$ satisfy the estimates
\begin{align*}
		\norm{\soluH}{\PVSHM{2}(\Gp)}&\leq C|k|^{-1}\left(\norm{\solf}{\PVSHM{1}(\Gp)}+|k|\norm{\solgi}{\VSHM{1/2}_T(\Gamma)}\right), \\
		\norm{\psi_{\solf}}{\PSHM{2}(\Gp)}&\leq C\norm{\diverg\solf}{\VSL(\Omega)}, \\
		\norm{\psi_{\solgi}}{\PSHM{2}(\Gp)}&\leq C\norm{\solgi}{\VSHM{1/2}_T(\Gamma)},
\end{align*}
where $C>0$ depends only on $\Gp$, $\mu^{-1}$, $\varepsilon$ and $\zeta$. 

Moreover, with $\rho(k)$ being the quantity from \eqref{rhodef} there holds
\begin{align*}
		\forall\ell\in\N_0:\ \norm{\soluA}{\PVSHM{\ell}(\Gp)}\leq \omega\left[\rho(k)+|k|^{-1}\right] \left(\norm{\solf}{\PVSHM{1}(\Gp)}+|k|\norm{\solgi}{\VSHM{1/2}_T(\Gamma)}\right)M^{\ell}(\ell+|k|)^{\ell},
\end{align*}
where the constants $\omega\geq 0, M>0$ depend only on $\Gp$, $\mu^{-1}$, $\varepsilon$ and $\zeta$.
In particular, $\soluA$ is piecewise analytic.
\end{theorem}

Our second main result considers a conforming Galerkin discretization of \eqref{Maxwellorig} based on N\'{e}d\'{e}lec-type I elements of degree $p$: Based on a (curvilinear) triangulation $\T$ with mesh width $h$ we use the discrete \ned-space $\Xh:=\nedelec(\T)$, see Section~\ref{discretizationsec} for the details. The discrete version of \eqref{weakformulation} then reads as follows: Find $\solu_N\in\Xh$ such that
\begin{align}\label{discreteproblem}
		\forall\solv_N\in\Xh:\ \bfa(\solu_N,\solv_N):= \SCP{\solf}{\solv_N}{\VSL(\Omega)}+\SCP{\solgi}{\Pi_T\solv_{N}}{\VSL_T(\Gamma)}.
\end{align}

Our second main result then states the following: The discrete problem \eqref{discreteproblem} has a unique solution $\solu_N$ and $\solu_N$ is quasi-optimal in $\Xh$ as long as there holds the $k$-explicit scale resolution condition that a) $|k|h/p$ is sufficiently small and b) $p/\log |k|$ is bounded from below.

\begin{theorem}\label{Mainresult2}
		Let $\Gp=\geom$ be an $\ana$-partition and let $\mu^{-1},\varepsilon\in\anatenspw{\Gp}$, $\zeta\in\anatens{\Gamma}$ and $\setu\subseteq\Coone$ be such that Assumption~\ref{assumption2} and Assumption ~\ref{assumption3} hold true. Furthermore, let $\T$ be a regular and shape regular (curvilinear) triangulation of $\Omega$ as described in Section~\ref{discretizationsec} below and suppose that $\T$ satisfies Assumption~\ref{Ttrafo}.

		Under these assumptions, consider a weak solution $\solu\in\HXI$ of Maxwell's equations \eqref{Maxwellorig} at wavenumber $k$ for some $k\in\setu$. 
		Then, there holds the following: For arbitrary $\eta\in (0,\cco^{-1})$ and $C_2>0$ there exists a constant $C_3>0$ depending only on $\eta,C_2,\Gp,\mu^{-1},\varepsilon,\zeta$, the parameters from Assumption~\ref{Ttrafo} and the constants $C_{\setu}>0$, \revision{$\theta\geq 0$} from \eqref{rhoalgebra} such that whenever we have
		\begin{align*}
				1+C_2\log |k|\leq p \quad {\rm and}\quad \frac{|k|h}{p}\leq C_3,
		\end{align*} 
		the corresponding discrete problem \eqref{discreteproblem} has a unique solution $\solu_N\in\Xh$, and
		\begin{align*}
				\norm{\solu-\solu_N}{\HXIK}\leq C\cco\frac{1+\eta}{1-\cco\eta}\ \inf_{\solw_N\in\Xh}\norm{\solu-\solw_N}{\HXIK}
		\end{align*}
		for a constant $C>0$ depending only on $\mu^{-1}, \varepsilon$ and $\zeta$.
\end{theorem}

		\begin{remark}
				Theorem~\ref{Mainresult2} relates the finite element error $\norm{\solu-\solu_N}{\HXIK}$ to the best approximation error $$\inf_{\solw_N\in\Xh}\norm{\solu-\solw_N}{\HXIK}.$$ As in \cite[Cor.~9.8]{MaxwellImpedanceMelenk} we may use the decomposition $\solu = \soluH+\soluA+k^{-2}\nabla\psi_{\solf}+k^{-1}\nabla\psi_{\solgi}$ from Theorem~\ref{Mainresult1} to estimate the best approximation error explicitly in $h$, $p$ and $k$. That is, in addition to providing quasi-optimality of the finite element error, Theorem~\ref{Mainresult1} and Theorem~\ref{Mainresult2} can be used as a starting point to provide wavenumber-explicit $hp$-FEM convergence estimates for Maxwell's equations \eqref{Maxwellorig}.
		\end{remark}


	
\section{Frequency splitting operators}\label{freqsec}

In this section we combine the idea of frequency splitting operators from, e.g., \cite{Esterhazy2012, MelenkHelmholtz1, HelmholtzMelenkSauter, MaxwellTransparentMelenk} with the regularity shift results from \cite{MaxwellMyself} to decompose solenoidal and piecewise regular vector fields into its high- and low-frequency components. 
These types of decompositions are the foundation of the regularity by decomposition technique from, e.g., \cite{Esterhazy2012, MelenkHelmholtz1, HelmholtzMelenkSauter,MaxwellTransparentMelenk} which will lead us to an elegant proof of Theorem~\ref{Mainresult1}.

\subsection{Frequency splitting operators in piecewise smooth media}

%
%

We continue by defining high- and low-frequency splitting operators acting on a smooth domain $\Omega$ consisting of $n\geq 1$ smooth subdomains. As a first step we consider the case $n=1$, and following e.g. \cite{MaxwellImpedanceMelenk} we make the following definition:

\begin{definition}\label{splittingglob}
		Let a wavenumber $k\in\Coone$ be given and let $\oldlambda>1$ be arbitrary. We define the full-space frequency splitting operators $\Higlob:\VSL(\R^3)\rightarrow\VSL(\R^3)$ and $\Lowglob:\VSL(\R^3)\rightarrow\VSL(\R^3)$ by 
	\begin{align*}
			\Higlob\solf := \F^{-1}\left((1-\chi_{\oldlambda |k|})\F[\solf]\right)\quad\rm{and}\ \ \Lowglob\solf := \F^{-1}\left(\chi_{\oldlambda |k|}\F[\solf]\right),
	\end{align*}
	where $\chi_{\oldlambda |k|}$ is the characteristic function of the open ball with radius $\oldlambda |k|$ centered at the origin, and the Fourier transformation $\F$ and its inverse $\F^{-1}$ are applied component-wise.
\end{definition}

By definition, $\Higlob$ and $\Lowglob$ are cutoff-operators \revtwo{in} Fourier space; since it is a priori not clear where the cutoff should take place we include the auxiliary parameter $\oldlambda>1$ which we will \revtwo{select} later according to our needs.

We notice that $\Lowglob\solf$ \revtwo{is band limited,} \revision{hence by the Paley-Wiener theorem} $\Lowglob\solf$ is analytic for all $\solf\in\VSL(\R^3)$. Moreover, $\Higlob$ is $\VSHM{\ell}(\R^3)$-stable and for $\solf\in\VSHM{\ell}(\R^3)$ one can estimate the growth of the derivatives of $\solf$ up to order $\ell$ explicitly in $\oldlambda$ and $k$, cf. \cite{HelmholtzMelenkSauter, MaxwellImpedanceMelenk}.

When applying these splitting operators to a piecewise regular vector field $\solf$ defined on a bounded domain $\Omega$ we have two issues: First, $\Higlob$ and $\Lowglob$ are only defined for vector fields on the full space, and second, it is not so clear whether $\Higlob\solv$ is again piecewise regular. Subsequently, we construct operators $\Higeom$ and $\Lowgeom$ which do not suffer from these issues at the cost of being defined only for solenoidal vector fields.

\medskip

In the proof of the subsequent Lemma~\ref{freqoplemma} we use the notions of jumps and normal jumps of functions and vector fields across subdomain interfaces. Let $\Gp = \geom$ be an $\ana$-partition and let $\interf_{i,j}:=\partial\Gp_i\cap\partial\Gp_j$ be a subdomain interface. Without loss of generality we assume $i<j$, and with $\soln_i$ denoting the outer unit normal of $\Gp_i$, we set
\begin{align}\label{jumpoperators}
		\jump{\varphi}\vert_{\interf_{i,j}} := (\varphi\vert_{\Gp_i})\vert_{\interf_{i,j}}-(\varphi\vert_{\Gp_j})\vert_{\interf_{i,j}} \quad{\rm and}\quad \njump{\solv}\vert_{\interf_{i,j}} := \left[(\solv\vert_{\Gp_i})\vert_{\interf_{i,j}}-(\solv\vert_{\Gp_j})\vert_{\interf_{i,j}}\right]\cdot\soln_i
\end{align}
for all piecewise sufficiently regular functions $\varphi$ and vector fields $\solv$. 
By Definition~\ref{partitiondef}, the union of all subdomain interfaces $\interf$ is the disjoint union of subdomain interfaces, hence \eqref{jumpoperators} induces jump operators $\jump{\cdot}$ and $\njump{\cdot}$ on $\interf$. 

\revision{In addition, for a bounded Lipschitz domain $\Omega\subseteq\R^3$ we define the space 
\begin{align*}
		\Hdivzero:=\left\{\solv\in\VSL(\Omega)\ |\ \diverg\solv = 0\ {\rm in}\ \Omega\right\}
\end{align*}
of all solenoidal fields on $\Omega$.}

%
%
\begin{lemma}\label{freqoplemma}
		Let $\Gp=\geom$ be an $\ana$-partition in the sense of Definition~\ref{partitiondef} and assume that a wavenumber $k\in\Coone$ is given. Let $\oldlambda>1$ be arbitrary. Under these assumptions, there exists an operator $\Higeom:\Hdivzero\rightarrow\Hdivzero$ with the following properties: For all $\ell\in\N_0$ and $0\leq \ell'\leq \ell$ and all $\solf\in\PVSHM{\ell}(\Gp)\cap\Hdivzero$ there holds
	\begin{align*}
		\norm{\Higeom\solf}{\PVSHM{\ell'}(\Gp)}\leq C_{\ell',\ell}(\oldlambda|k|)^{\ell'-\ell}\norm{\solf}{\PVSHM{\ell}(\Gp)}
	\end{align*}	
for a constant $C_{\ell',\ell}>0$ independent of $k, \solf$ and $\oldlambda$. 

Moreover, there exists an operator $\Lowgeom:\Hdivzero\rightarrow\Hdivzero$ which for all $\solf\in\Hdivzero$ satisfies $\Higeom\solf+\Lowgeom\solf = \solf$ and
	\begin{align*}
			\forall\ell\in\N_0:\ \norm{\Lowgeom\solf}{\PVSHM{\ell}(\Gp)}\leq \omega M^{\ell}\oldlambda^{\ell}\left(\ell+|k|\right)^{\ell}\norm{\solf}{\VSL(\Omega)}
	\end{align*}
	for constants that are $\omega>0,M>0$ are independent of $k$, $\oldlambda$, $\solf$ and $\ell$. In particular, there holds that $\Lowgeom\solf\in\anaomegavecpw{\Gp}{\omega\norm{\solf}{\VSL(\Omega)}}{M\oldlambda}$. 
\end{lemma}

\begin{fatproof}
		We divide the proof into two steps. In the first step we construct frequency splitting operators on every single subdomain and in the second step we conclude.

		\textbf{Step 1:} \revision{We consider a single subdomain $\Gp_i$ for some $i\in\{0,\ldots,n\}$. We claim that there exist operators $\Higeomdom{\Gp_i},\Lowgeomdom{\Gp_i}:\Hdivzerodom{\Gp_i}\rightarrow\Hdivzerodom{\Gp_i}$ which for all $\solf\in\Hdivzerodom{\Gp_i}$ satisfy $\Higeomdom{\Gp_i}\solf+\Lowgeomdom{\Gp_i}\solf = \solf$ and
		\begin{equation}\label{freqsplitgpi}
		\begin{alignedat}{2}
				\forall 0\leq \ell'\leq\ell:\ &\norm{\Higeomdom{\Gp_i}\solf}{\VSHM{\ell'}(\Gp_i)}&&\leq C_{\ell',\ell}(\oldlambda|k|)^{\ell'-\ell}\norm{\solf}{\VSHM{\ell}(\Gp_i)}, \\
				\forall\ell\in\N_0:\ &\norm{\Lowgeomdom{\Gp_i}\solf}{\VSHM{\ell}(\Gp_i)}&&\leq \omega M^{\ell}\oldlambda^{\ell}\left(\ell+|k|\right)^{\ell}\norm{\solf}{\VSL(\Gp_i)}.
		\end{alignedat}
\end{equation}
		 We mimic the construction in \cite[Sec.~6]{MaxwellImpedanceMelenk}: According to e.g. \cite[Lem.~2.6]{MaxwellImpedanceMelenk} or \cite[Cor.~4.8]{CurlInverse} there exist operators $\Roperator{\Gp_i}, \Koperator{\Gp_i}$ which for all integers $m\in\Z$ have the mapping properties $\Roperator{\Gp_i}:\VSHM{-m}(\Gp_i)\rightarrow\VSHM{1-m}(\Gp_i)$ and $\Koperator{\Gp_i}:\VSHM{m}(\Gp_i)\rightarrow\left(\CM{\infty}(\overline{\Gp_i})\right)^3$ and satisfy
		\begin{align*}
				\forall\solf\in\Hdivzerodom{\Gp_i}:\ \solf = \curl\Roperator{\Gp_i}\solf+\Koperator{\Gp_i}\solf.
		\end{align*}
		Using Stein's extension operator $\stein{\Gp_i}$ from \cite{Stein} we define
		\begin{align}\label{tildesplitting}
\begin{split}
				\Higeomdomt{\Gp_i}\solf&:=\left(\curl\Higlob\stein{\Gp_i}\Roperator{\Gp_i}\solf\right)\vert_{\Gp_i}+\left(\Higlob\stein{\Gp_i}\Koperator{\Gp_i}\solf\right)\vert_{\Gp_i}, \\
				\Lowgeomdomt{\Gp_i}\solf&:=\left(\curl\Lowglob\stein{\Gp_i}\Roperator{\Gp_i}\solf\right)\vert_{\Gp_i}+\left(\Lowglob\stein{\Gp_i}\Koperator{\Gp_i}\solf\right)\vert_{\Gp_i}.
\end{split}
		\end{align}
		By construction, $\Higeomdomt{\Gp_i}\solf+\Lowgeomdomt{\Gp_i}\solf = \solf$ for all $\solf\in\Hdivzero$ and with similar techniques as in \cite[Lem.~6.5]{MaxwellImpedanceMelenk} we get that\footnote{See \cite[Ch.~9]{MyThesis} for further details.}
		\begin{equation}\label{freqtildestimates}
			\begin{alignedat}{2}
					\forall 0\leq \ell'\leq\ell:\ &\norm{\Higeomdomt{\Gp_i}\solf}{\VSHM{\ell'}(\Gp_i)}&&\leq C_{\ell',\ell}(\oldlambda|k|)^{\ell'-\ell}\norm{\solf}{\VSHM{\ell}(\Gp_i)}, \\
					\forall\ell\in\N_0:\ &\norm{\Lowgeomdomt{\Gp_i}\solf}{\VSHM{\ell}(\Gp_i)}&&\leq \omega M^{\ell}\oldlambda^{\ell} (\ell+|k|)^{\ell}\norm{\solf}{\VSL(\Gp_i)}
					\end{alignedat}
		\end{equation}
for constants $\omega, M>0$ independent of $\lambda$, $k$ and $\ell$.}
		
\revision{			However, although we suppose $\diverg\solf=0$ it is not guaranteed that $\diverg\Higeomdomt{\Gp_i}\solf=0$ or $\diverg\Lowgeomdomt{\Gp_i}\solf=0$. Therefore we slightly modify the definitions \eqref{tildesplitting} and set
		\begin{alignat*}{2}
				\Higeomdom{\Gp_i}\solf&:=\Higeomdomt{\Gp_i}\solf+\nabla\psi_{\solf,H,\Gp_i}, \\
				\Lowgeomdom{\Gp_i}\solf&:=\Lowgeomdomt{\Gp_i}\solf+\nabla\psi_{\solf,L,\Gp_i},
		\end{alignat*}
where $\psi_{\solf, H}$ and $\psi_{\solf,L}$ are the solutions of}

\revision{
\begin{minipage}{0.5\textwidth}
		\begin{alignat*}{2}
				-\diverg\nabla\psi_{\solf,H,\Gp_i} &= \diverg\Higeomdomt{\Gp_i}\solf &&\quad{\rm in}\quad \Gp_i,\\
				\psi_{\solf,H,\Gp_i}&=0 && \quad{\rm on}\quad\partial\Gp_i,
	\end{alignat*}
\end{minipage}
\begin{minipage}{0.5\textwidth}
		\begin{alignat*}{2}
				-\diverg\nabla\psi_{\solf,L,\Gp_i} &= \diverg\Lowgeomdomt{\Gp_i}\solf &&\quad{\rm in}\quad \Gp_i,\\
				\psi_{\solf,L,\Gp_i}&=0 && \quad{\rm on}\quad\partial\Gp_i.
	\end{alignat*}
\end{minipage}
}

\bigskip

\revision{By construction, $\Higeomdom{\Gp_i},\Lowgeomdom{\Gp_i}:\Hdivzerodom{\Gp_i}\rightarrow\Hdivzerodom{\Gp_i}$. Furthermore, due to $\Higeomdomt{\Gp_i}\solf+\Lowgeomdomt{\Gp_i}\solf=\solf$ we have we have $\nabla\psi_{\solf,H,\Gp_i}+\nabla\psi_{\solf,L,\Gp_i}=0$ for all $\solf\in\Hdivzerodom{\Gp_i}$. This implies $\Higeomdom{\Gp_i}\solf+\Lowgeomdom{\Gp_i}\solf = \solf$ for all $\solf\in\Hdivzerodom{\Gp_i}$. It remains to argue \eqref{freqsplitgpi}. In fact, there holds
		\begin{equation}\label{eq:tmpcomput}
		\begin{alignedat}{2}
		\forall 0\leq \ell'\leq \ell:\ &\norm{\nabla\psi_{\solf,H,\Gp_i}}{\VSHM{\ell'}(\Gp_i)}&&\leq C_{\ell'}\norm{\Higeomdomt{\Gp_i}\solf}{\VSHM{\ell'}(\Gp_i)}, \\
		\forall\ell\in\N_0:\ &\norm{\nabla\psi_{\solf,L,\Gp_i}}{\VSHM{\ell}(\Gp_i)}&&\leq \omega(\oldlambda M)^{\ell}(\ell+|k|)^{\ell}\norm{\solf}{\VSL(\Gp_i)}
\end{alignedat}
\end{equation}
for constants $\omega,M>0$ independent of $\oldlambda$, $k$ and $\ell$. For $\nabla\psi_{\solf,H,\Gp_i}$ this follows from \cite[Thm.~2.3.2, Thm.~3.4.1]{GrandLivre} and the Lax-Milgram lemma. }

\revtwo{For $\nabla\psi_{\solf,L,\Gp_i}$ we use \cite[Thm.~3.4.2]{GrandLivre} and arguments from \cite[Thm.~2.7.1]{GrandLivre} to obtain
\begin{align*}
		\forall\ell\in\N:\ \max_{|\alpha|=\ell+1}\norm{\D^{\alpha}(\psi_{\solf,L,\Gp_i})}{\SL(\Gp_i)}\leq A^{\ell+2}\left(\norm{\psi_{\solf,L,\Gp_i}}{\VSHM{1}(\Gp_i)}+\sum_{m=0}^{\ell-1}\frac{(\ell+1)!}{m!}\norm{\Lowgeomdomt{\Gp_i}\solf}{\VSHM{m}(\Gp_i)}\right),
\end{align*}
and by exploiting $\ell+1\leq 2^{\ell}$ this becomes
\begin{align*}
		\forall\ell\in\N:\ \max_{|\alpha|=\ell+1}\norm{\D^{\alpha}(\psi_{\solf,L,\Gp_i})}{\SL(\Gp_i)}\leq A^{\ell+2}\left(\norm{\psi_{\solf,L,\Gp_i}}{\VSHM{1}(\Gp_i)}+\sum_{m=0}^{\ell-1}\frac{(\ell)!}{m!}\norm{\Lowgeomdomt{\Gp_i}\solf}{\VSHM{m}(\Gp_i)}\right)
\end{align*}
for a modified constant $A\geq 1$ that is independent of $k$, $\ell$ and $\oldlambda$.
By Stirling's formula and \cite[Lem.~A.1]{MaxwellMyself2} we obtain
\begin{align*}
		\forall\ell\in\N:\ \max_{|\alpha|=\ell}\norm{\D^{\alpha}(\psi_{\solf,L,\Gp_i})}{\SL(\Gp_i)}\leq A^{\ell+2}\left(\norm{\psi_{\solf,L,\Gp_i}}{\VSHM{1}(\Gp_i)}+\sum_{m=0}^{\ell-1}\frac{(\ell+|k|)^{\ell}}{(m+|k|)^m}\norm{\Lowgeomdomt{\Gp_i}\solf}{\VSHM{m}(\Gp_i)}\right)
\end{align*}
for a constant $A\geq 1$ that is again independent of $k$, $\oldlambda$ and $\ell$.}
\revtwo{As a consequence of the Lax-Milgram lemma we have $$\norm{\psi_{\solf,L,\Gp_i}}{\SL(\Gp_i)}\leq C\norm{\solf}{\VSL(\Gp_i)},$$ and in combination with \eqref{freqtildestimates} this leads to 
\begin{align*}
		\forall\ell\in\N:\ \max_{|\alpha|=\ell+1}\norm{\D^{\alpha}(\psi_{\solf,L,\Gp_i})}{\SL(\Gp_i)}&\leq \omega A^{\ell+2}(\ell+|k|)^{\ell}\norm{\solf}{\VSL(\Gp_i)}\sum_{m=0}^{\ell-1}(M\oldlambda)^m(m+|k|)^m \\
																								   &\leq \omega A^{\ell+2}(M\oldlambda)^{\ell}(\ell+|k|)^{\ell} \norm{\solf}{\VSL(\Gp_i)}\underbrace{\sum_{m=0}^{\ell-1}1}_{\leq 2^{\ell}}
\end{align*}
where we exploited $M\oldlambda\geq1$ for the second inequality. In total, we get
\begin{align*}
		\forall\ell\in\N:\ \max_{|\alpha|=\ell}\norm{\D^{\alpha}(\psi_{\solf,L,\Gp_i})}{\SL(\Gp_i)}
		\leq \omega (M\oldlambda)^{\ell}(\ell+|k|)^{\ell}\norm{\solf}{\VSL(\Gp_i)} 
\end{align*}
for constants $\omega>0$, $M\geq 1$ independent of $k$, $\oldlambda$ and $\ell$.}
Finally, we have
\begin{align*}
		\norm{\nabla\psi_{\solf,L,\Gp_i}}{\VSHM{\ell}(\Gp_i)}^2 = \sum_{m=0}^{\ell+1}\max_{|\alpha|=m}\norm{\D^{\alpha}\psi_{\solf,L,\Gp_i}}{\SL(\Gp_i)}^2 
														 &\leq \omega^2\norm{\solf}{\VSL(\Gp_i)}^2\sum_{m=0}^{\ell+1}(M \oldlambda)^{2\ell}(m+|k|)^{2m} \\
														 &\leq \omega^2\norm{\solf}{\VSL(\Gp_i)}^2(M\oldlambda)^{2\ell}(\ell+|k|)^{2\ell}\underbrace{\sum_{m=0}^{\ell+1} 1}_{\leq 2^{\ell+1}},
\end{align*}
that is,
\begin{align*}
\forall\ell\in\N_0:\ \norm{\nabla\psi_{\solf,L,\Gp_i}}{\VSHM{\ell}(\Gp_i)}\leq \omega(M\oldlambda)^{\ell}(\ell+|k|)^{\ell}\norm{\solf}{\VSL(\Gp_i)}
\end{align*}
for constants $\omega>0$, $M\geq 1$ independent of $k$, $\oldlambda$ and $\ell$. This concludes the proof of \eqref{eq:tmpcomput} and together with \eqref{freqtildestimates} we conclude \eqref{freqsplitgpi}.

\bigskip

\textbf{Step 2:} After having defined frequency splitting operators $\Higeomdom{\Gp_i}, \Lowgeomdom{\Gp_i}$ on a single subdomain $\Gp_i$, we construct frequency splitting operators on an $\ana$-partition $\Gp = \geom$. We do this by first defining operators $\Higeomstar$ and $\Lowgeomstar$ piecewise on every subdomain by $\Higeomstar\solf\vert_{\Gp_i}:=\Higeomdom{\Gp_i}(\solf\vert_{\Gp_i})$ and $\Lowgeomstar\solf\vert_{\Gp_i}:=\Lowgeomdom{\Gp_i}(\solf\vert_{\Gp_i})$, respectively. Although these operators inherit the nice mapping properties from \eqref{freqsplitgpi}, it is not guaranteed that the resulting vector field is in $\Hdiv$. 

		\medskip

The remedy is to modify the previous definition of $\Higeomstar$ and $\Lowgeomstar$ by adding a suitable gradient field which lifts possible normal jumps.
	To that end we define splitting operators $\PsiH:\Hdivzero\rightarrow\VSL(\Omega)$ and $\PsiL:\Hdivzero\rightarrow\VSL(\Omega)$ by $\PsiH\solf:=\nabla\varphi_H$ and $\PsiL\solf:=\nabla\varphi_L$, where $\varphi_H$ and $\varphi_L$ are solutions of the transmission problems

\begin{minipage}{0.5\textwidth}
		\begin{alignat*}{2}
				\Delta(\varphi_H\vert_{\Gp_i}) &= 0 &&\quad{\rm in}\quad \Gp_1,...,\Gp_n,\\
				\jump{\varphi_H}&=0 && \quad {\rm on}\quad \interf,\\
				\njump{\nabla\varphi_H} &= \njump{\Higeomstar\solf}&& \quad{\rm on}\quad \interf,\\
				\varphi_H&=0 && \quad{\rm on}\quad\Gamma,
	\end{alignat*}
\end{minipage}
\begin{minipage}{0.5\textwidth}
		\begin{alignat*}{2}
			\Delta(\varphi_L\vert_{\Gp_i}) &= 0 &&\quad{\rm in}\quad \Gp_1,...,\Gp_n,\\
			\jump{\varphi_L}&=0 &&\quad{\rm on}\quad \interf,\\
	\njump{\nabla\varphi_L} &= \njump{\Lowgeomstar\solf} &&\quad{\rm on}\quad \interf,\\
	\varphi_L&=0 &&\quad{\rm on}\quad \Gamma.
	\end{alignat*}
\end{minipage}

\medskip

We notice that existence and uniqueness of $\varphi_H$ and $\varphi_L$ follow from the Lax-Milgram~lemma, and \revision{for all $\solf\in\Hdivzero$ we have $\njump{\Higeomstar\solf}+\njump{\Lowgeomstar\solf}=0$ on $\interf$, hence $\PsiH\solf+\PsiL\solf = 0$ for all $\solf\in\Hdivzero$.} Furthermore, we claim that the operators $\PsiH$ and $\PsiL$ have nice mapping properties: Indeed, for $0\leq \ell'\leq \ell$ there holds
\begin{align*}
		\norm{\PsiH\solf}{\PVSHM{\ell'}(\Omega)}\leq C_{\ell'}\norm{\Higeomstar\solf}{\PVSHM{\ell'}(\Omega)}
		\end{align*}
		for a constant $C_{\ell'}>0$ depending only on $\Gp$ and $\ell'$. For $\ell'=0$ this follows from the Lax-Milgram lemma and a trace inequality, and for $\ell'>0$ this follows from \cite[Thm.~5.3.8]{GrandLivre} and again a trace inequality. 

		In a similar way, \cite[Thm.~5.3.8, (3.32)]{GrandLivre} and a \revtwo{similar computation as in the treatment of $\psi_{\solf,L,\Gp_i}$ above} show 
        \begin{align*}
				\forall\ell\in\N_0:\ \norm{\PsiL\solf}{\PVSHM{\ell}(\Omega)}\leq \omega (\oldlambda M)^{\ell}\left(\ell+|k|\right)^{\ell}\norm{\Lowgeomstar\solf}{\VSL(\Omega)}
		\end{align*}
		for constants $\omega,M>0$ independent of $\oldlambda$, $k$ and $\ell$.

\medskip

Finally, for any $\solf\in\Hdivzero$ we define 
\begin{align*}
		\Higeom\solf:=\Higeomstar\solf-\PsiH\solf\quad {\rm and}\quad \Lowgeom\solf:=\Lowgeomstar\solf-\PsiL\solf.
\end{align*}
By construction, $\Higeom$ and $\Lowgeom$ map $\Hdivzero$ to itself and satisfy $\Higeom\solf+\Lowgeom\solf = \solf$ for all $\solf\in\Hdivzero$. Furthermore, $\Higeom$ and $\Lowgeom$ inherit the sought mapping properties from \eqref{freqsplitgpi}. This concludes the proof.

\end{fatproof}

\subsection{Frequency splitting operators on analytic surfaces}

In addition to frequency splitting operators on $\ana$-partitions we will also need splitting operators that decompose boundary traces of vector fields into its high- and low-frequency components. Thankfully, the situation on the boundary is easier than the case of splitting operators on $\ana$-partitions, since we do not need to account for interface conditions between subdomains.

The following proposition is due to \cite[Sec.~6.3]{MaxwellImpedanceMelenk}; it asserts the existence of suitable high- and low-frequency splittings of tangent fields on analytic surfaces. 
\begin{proposition}\label{freqopbdrprop}
		Let $\Gp=\geom$ be an $\ana$-partition and let $\Gamma$ denote the (by definition analytic) boundary of $\Omega$. Furthermore, let $k\in\Coone$ be a given wavenumber and let $\oldlambda > 1$ be arbitrary.
		Then, there exists an operator $\Higamma$ with the following mapping properties: For any $0\leq s'\leq s$ there exist constants $C_{s',s}>0$ independent of $\oldlambda, \solg$ and $k$ such that for all $\solg\in\VSHM{s}_T(\Gamma)$
			\begin{align*}
		\norm{\Higamma\solg}{\VSHM{s'}_T(\Gamma)}&\leq C_{s',s}(\oldlambda|k|)^{s'-s}\norm{\solg}{\VSHM{s}_T(\Gamma)}, \\
		\norm{\diverg_{\Gamma}\Higamma\solg}{\SHM{s'-1}(\Gamma)}&\leq C_{s',s}(\oldlambda|k|)^{s'-s}\norm{\diverg_{\Gamma}\solg}{\SHM{s-1}(\Gamma)}, \\
		\norm{\curl_{\Gamma}\Higamma\solg}{\SHM{s'-1}(\Gamma)}&\leq C_{s',s}(\oldlambda|k|)^{s'-s}\norm{\curl_{\Gamma}\solg}{\SHM{s-1}(\Gamma)}.
	\end{align*}

	Moreover, there exists an operator $\Logamma$ such that for all $\solg\in\VSL_T(\Gamma)$ there holds $\Higamma\solg+\Logamma\solg=\solg$ and 
	\begin{align*}
			\Logamma\solg\in\anagammavecspecial{\Gamma}{\omega\norm{\solg}{\VSL_T(\Gamma)}}{M}, 
	\end{align*}
	where $\omega>0$ is independent of $k, \oldlambda$ and $\solg$, and $M>0$ is independent of $k$ and $\solg$ (but may depend on $\oldlambda$).
	
\end{proposition}

\section{The auxiliary Maxwell problem}\label{auxmax}

One main difficulty in the study of the Maxwell problem \eqref{Maxwellorig} is that wavenumber-explicit a priori estimates on the solution $\solu$ are hard to establish. Even if we assumed that Assumption~\ref{assumption3} holds, we could estimate the amplifying factor $\rho(k)$ from \eqref{rhodef} only in terms of $|k|^{\theta}$ for some generally unknown $\theta\geq 0$, which is too weak for our purposes.

We circumvent this issue by considering an auxiliary Maxwell problem which fits into the Lax-Milgram setting, an idea which goes back to \cite{HelmholtzMelenkSauter} for the Helmholtz equation and \cite{MaxwellTransparentMelenk, MaxwellImpedanceMelenk} for Maxwell's equations. 

\medskip

Let $\Gp = \geom$ be an $\ana$-partition and suppose that the coefficients $\mu^{-1},\varepsilon\in\anatenspw{\Gp}$ satisfy \eqref{coercivemu}-\eqref{coerciveeps} and that $\zeta\in\anagamma$ satisfies \eqref{coercivezeta} and \eqref{zetaproperty}. 
Due to these assumptions, we have that there exists a constant $c>0$ and unimodular numbers $\alpha_{\varepsilon}$ and $\alpha_{\zeta}$ such that 
\begin{align*}
		\realpart\SCP{\mu^{-1}\solz}{\solz}{}\geq c\norm{\solz}{}^2, \quad \realpart\SCP{\alpha_{\varepsilon}\varepsilon\solz}{\solz}{}\geq c\norm{\solz}{}^2\quad{\rm and}\quad \realpart\SCP{\alpha_{\zeta}\zeta\solz}{\solz}{}\geq c\norm{\solz}{}^2
\end{align*}
for all $\solz\in\Co^{3}$ uniformly in $\Omega$ or on the boundary $\Gamma$, respectively.
Depending on $\alpha_{\varepsilon}$ and $\alpha_{\zeta}$ and an arbitrary but fixed wavenumber $k\in\Coone$ we set
\begin{align*}
		\sigma_{\varepsilon}:=-\frac{\alpha_{\varepsilon}|k|^2}{k^2}\quad{\rm and}\quad \sigma_{\zeta}:=\frac{i\alpha_{\zeta}|k|}{k},
\end{align*}
\revtwo{and} notice that $\sigma_{\varepsilon}$ and $\sigma_{\zeta}$ are again unimodular.

Based on $\sigma_{\varepsilon}$ and $\sigma_{\zeta}$ we define an auxiliary Maxwell problem, which reads as follows: Find a vector field $\soluaux$ such that 

\begin{equation}\label{Maxwellauximp}
		\begin{alignedat}{2}
				\curl\mu^{-1}\curl\soluaux-\sigma_{\varepsilon}k^2\varepsilon\soluaux &= \solf \quad &&{\rm in}\quad \Omega, \\
				\left(\mu^{-1}\curl\soluaux\right)\times\soln-ik\sigma_{\zeta}\zeta\soluaux_T &= \solgi \quad &&{\rm on}\quad \Gamma.
\end{alignedat}
\end{equation}
With the space $\HXI$ from \eqref{energyspace}, the weak formulation of \eqref{Maxwellauximp} reads as follows: Find $\soluaux\in\HXI$ such that
\begin{align*}
		\SCP{\mu^{-1}\curl \soluaux}{\curl\solv}{\VSL(\Omega)}-k^2\sigma_{\varepsilon}\SCP{\varepsilon\soluaux}{\solv}{\VSL(\Omega)}-ik\sigma_{\zeta}\SCP{\zeta\soluaux_T}{\solv_T}{\VSL_T(\Gamma)}
					= \SCP{\solf}{\solv}{\VSL(\Omega)}+\SCP{\solgi}{\solv_T}{\VSL_T(\Gamma)}
\end{align*}
for all $\solv\in\HXI$. Existence and uniqueness of a weak solution $\soluaux$ of \eqref{Maxwellauximp} follows from Lemma~\ref{laxmilgram} below.

\subsection{Well-posedness of the auxiliary Maxwell problem}

The main difference between the auxiliary problem \eqref{Maxwellauximp} and \revision{the original problem} \eqref{Maxwellorig} is that the auxiliary problem fits into the Lax-Milgram setting. That is, we have existence and uniqueness of weak solutions together with $k$-explicit dependencies on the given data, as the following lemma shows:

\begin{lemma}\label{laxmilgram}
		Let $\Gp$ be an $\ana$-partition and suppose that $k\in\Coone$ is given. Furthermore, suppose that the coefficients $\mu^{-1},\varepsilon\in\anatenspw{\Gp}$ satisfy \eqref{coercivemu}-\eqref{coerciveeps} and that $\zeta\in\anagamma$ satisfies \eqref{coercivezeta} and \eqref{zetaproperty}. 

		\smallskip

		Under these assumptions, let $\solf\in\VSL(\Omega)$ and boundary data $\solgi\in\VSL_T(\Gamma)$ be given. Then, the corresponding auxiliary problem \eqref{Maxwellauximp} has a unique weak solutions $\soluaux$, which satisfies 
		\begin{align}\label{soluauxestimate}
				\norm{\soluaux}{\VSL(\Omega)}+|k|^{-1}\norm{\curl\soluaux}{\VSL(\Omega)}+|k|^{-1/2}\norm{\soluaux_T}{\VSL_T(\Gamma)}\leq C\left(|k|^{-2}\norm{\solf}{\VSL(\Omega)}+|k|^{-3/2}\norm{\solgi}{\VSL_T(\Gamma)}\right),
		\end{align}
		where the constant $C>0$ depends only on $\Omega$, $\mu^{-1}$, $\varepsilon$ and $\zeta$.
\end{lemma}

\begin{fatproof}
				By construction of $\sigma_{\varepsilon}$ and $\sigma_{\zeta}$, the sesquilinear form 
				\begin{align*}
						A_k^{\sigma}(\solv,\solw) := \SCP{\mu^{-1}\curl \solv}{\curl\solw}{\VSL(\Omega)}-k^2\sigma_{\varepsilon}\SCP{\varepsilon\solv}{\solw}{\VSL(\Omega)}-ik\sigma_{\zeta}\SCP{\zeta\solv_T}{\solw_T}{\VSL_T(\Gamma)}
				\end{align*}
				satisfies 
				\begin{align*}
						\forall\solv\in\HXI:\ \realpart A_k^{\sigma}(\solv,\solv)\geq c\left(\norm{\curl\solv}{\VSL(\Omega)}^2+|k|^2\norm{\solv}{\VSL(\Omega)}^2+|k|\norm{\solv_T}{\VSL_T(\Gamma)}^2\right),
				\end{align*}
				where the constant $c>0$ depends only on $\Omega$, $\mu^{-1}$, $\varepsilon$ and $\zeta$. Hence, the statement follows from the Lax-Milgram lemma.

\end{fatproof}

Under additional regularity assumptions on the given data higher regularity estimates on $\solu$ are possible. The subsequent result establishes piecewise $\VSHM{1}$-regularity of $\solu$ and $\curl\solu$ provided that the given data is a little bit more regular. For the sake of simplicity we restricted ourselves to the case of solenoidal $\solf$.

\begin{lemma}\label{regshift1}
		Let $\Gp$ be an $\ana$-partition and suppose that $k\in\Coone$ is given. 
		Furthermore, assume that the coefficients $\mu^{-1}, \varepsilon\in\anatenspw{\Gamma}$ satisfy \eqref{coercivemu}-\eqref{coerciveeps} and that $\zeta\in\anagamma$ satisfies \eqref{coercivezeta} and \eqref{zetaproperty}. 

		In addition, suppose that $\solf\in\VSL(\Omega)$ satisfies $\diverg\solf = 0$ and that $\solgi\in\VSHM{1/2}_T(\Gamma)$.
Then, the weak solution $\soluaux$ of the auxiliary problem \eqref{Maxwellauximp} satisfies $\soluaux\in\PVSHM{1}(\Gp)\cap\PVHcurl{1}$ together with the estimate
\begin{align*}
		\norm{\soluaux}{\PVSHM{1}(\Gp)}+|k|^{-1}\norm{\curl\soluaux}{\PVSHM{1}(\Gp)}\leq C|k|^{-1}\left(\norm{\solf}{\VSL(\Omega)}+|k|^{1/2}\norm{\solgi}{\VSL_T(\Gamma)}+ \norm{\solgi}{\VSHM{1/2}_T(\Gamma)}\right),
\end{align*}
where the constant $C>0$ depends only on $\Gp$, $\mu^{-1}$, $\varepsilon$ and $\zeta$.

\end{lemma}

\begin{fatproof}
		Follows from Lemma~\ref{laxmilgram} and \cite[Thm.~2.10]{MaxwellMyself}.
\end{fatproof}

If $\solf$ and $\solgi$ are even more regular we may conclude $\PVSHM{2}$-regularity of $\soluaux$, as we shall see next:

\begin{lemma}\label{regshift2}
		Let $\Gp$ be an $\ana$-partition and assume that $k\in\Coone$ is given. Furthermore, suppose that the coefficients $\mu^{-1}, \varepsilon\in\anatenspw{\Gamma}$ satisfy \eqref{coercivemu}-\eqref{coerciveeps} and that $\zeta\in\anagamma$ satisfies \eqref{coercivezeta} and \eqref{zetaproperty}. 

		In addition, suppose that $\solf\in\VSL(\Omega)$ satisfies $\diverg\solf = 0$, as well as $\solgi\in\VSHM{1/2}_T(\Gamma)$ and $\solf\cdot\soln-\diverg_{\Gamma}\solgi\in\SHM{1/2}(\Gamma)$.
Then, the weak solution $\soluaux$ of the auxiliary problem \eqref{Maxwellauximp} satisfies $\soluaux\in\PVSHM{2}(\Gp)$ and there holds
\begin{align*}
				\norm{\soluaux}{\PVSHM{2}(\Gp)} \leq C\left(\norm{\solf}{\VSL(\Omega)}+|k|^{1/2}\norm{\solgi}{\VSL_T(\Gamma)}+\norm{\solgi}{\VSHM{1/2}_T(\Gamma)}+|k|^{-1}\norm{\solf\cdot\soln-\diverg_{\Gamma}\solgi}{\SHM{1/2}(\Gamma)}\right),
\end{align*}
where the constant $C>0$ depends only on $\Gp$, $\mu^{-1}$, $\varepsilon$ and $\zeta$.
\end{lemma}

\begin{fatproof}
		The proof follows the lines of \cite[Proof of Thm.~2.10, Step 1]{MaxwellMyself}. We recall that for $\solv\in\Hcurl$ there holds the identity 
		\begin{align}\label{surfaceidentities}
				\diverg_{\Gamma}\solv_t = \curl_{\Gamma}\solv_T = (\curl\solv)\vert_{\Gamma}\cdot\soln,
		\end{align}
		where $\soln$ denotes the outer unit normal to $\Gamma$. According to \cite[Thm.~2.7]{MaxwellMyself} we have
\begin{align*}
		\norm{\soluaux}{\PVSHM{2}(\Gp)}\leq C\left(\norm{\curl\soluaux}{\PVSHM{1}(\Gp)}+\norm{\diverg\varepsilon\soluaux}{\PSHM{1}(\Omega)}+\norm{\soluaux_T}{\VSHM{3/2}_T(\Gamma)}\right),
\end{align*}
and due to $\diverg\varepsilon\soluaux = -\sigma_{\varepsilon}^{-1}k^{-2}\diverg\solf = 0$ we arrive at
	\begin{align*}
		\norm{\soluaux}{\PVSHM{2}(\Gp)}\leq C\left(\norm{\curl\soluaux}{\PVSHM{1}(\Gp)}+\norm{\soluaux_T}{\VSHM{3/2}_T(\Gamma)}\right).
	\end{align*}
	Furthermore, according to \cite[Lem.~3.6]{MaxwellMyself} there holds
	\begin{align*}
			\norm{\soluaux_T}{\VSHM{3/2}_T(\Gamma)}\leq C\left(\norm{\diverg_{\Gamma}\zeta\soluaux_T}{\SHM{1/2}(\Gamma)}+\norm{\curl_{\Gamma}\soluaux_T}{\SHM{1/2}(\Gamma)}\right).
	\end{align*}
	According to \eqref{surfaceidentities} we have $\curl_{\Gamma}\soluaux_T = \curl\soluaux\cdot\soln$ on $\Gamma$, hence $\norm{\curl_{\Gamma}\soluaux_T}{\SHM{1/2}(\Gamma)}\leq C\norm{\curl\solu}{\PVSHM{1}(\Gp)}$. In addition, taking the surface divergence on the impedance boundary condition in \eqref{Maxwellauximp} yields
	\begin{align*}
			ik\sigma_{\zeta}\diverg_{\Gamma}\soluaux_T &= \diverg_{\Gamma}(\mu^{-1}\curl\soluaux\times\soln)-\diverg_{\Gamma}\solgi \\
												 &= \curl\mu^{-1}\curl\soluaux\cdot\soln-\diverg_{\Gamma}\solgi \\
												 &= \solf\cdot\soln+\sigma_{\varepsilon}k^2\varepsilon\soluaux\cdot\soln-\diverg_{\Gamma}\solgi,
	\end{align*}
	hence
\begin{align*}
		\norm{\soluaux_T}{\VSHM{3/2}_T(\Gamma)}\leq C\left(\norm{\curl\soluaux}{\PVSHM{1}(\Gp)}+|k|^{-1}\norm{\solf\cdot\soln-\diverg_{\Gamma}\solgi}{\SHM{1/2}(\Gamma)}+|k|\norm{\soluaux}{\PVSHM{1}(\Gp)}\right).
\end{align*}
In total we get 
\begin{align*}
		\norm{\soluaux}{\PVSHM{2}(\Gp)}\leq C\left(\norm{\curl\soluaux}{\PVSHM{1}(\Gp)}+|k|^{-1}\norm{\solf\cdot\soln-\diverg_{\Gamma}\solgi}{\SHM{1/2}(\Gamma)}+|k|\norm{\soluaux}{\PVSHM{1}(\Gp)}\right),
\end{align*}
and together with Lemma~\ref{regshift1} this concludes the proof. 
\end{fatproof}

\section{Proof of Theorem~\ref{Mainresult1}}\label{Mainres1proof}

Before presenting a proof of Theorem~\ref{Mainresult1} we show two auxiliary results.
The first of these two results is similar to \cite[Lem.~7.2]{MaxwellImpedanceMelenk}, it allows us to lift surface divergences on $\Gamma$ to gradient fields in $\Omega$. 

\begin{lemma}\label{liftg}
		Let $\Gp$ be an $\ana$-partition and suppose that the coefficient $\varepsilon\in\anatenspw{\Gp}$ satisfies \eqref{coerciveeps} and that $\zeta\in\anagamma$ satisfies \eqref{coercivezeta} and \eqref{zetaproperty}. Then, for all $\ell\in\N_0$ and every tangent field $\solg\in\VSHM{\ell+1/2}_T(\Gamma)$ there exists a function $\psi_{\solg}\in\PSHM{\ell+2}(\Gp)$ which satisfies
		\begin{equation*}
				\begin{alignedat}{2}
					\diverg\varepsilon\nabla\psi_{\solg} &=0 \quad &&{\rm in} \quad \Omega, \\
				\diverg_{\Gamma}\zeta\nabla_{\Gamma}\psi_{\solg} & = -i\diverg_{\Gamma}\solg \quad &&{\rm on} \quad \Gamma,
				\end{alignedat}
					\end{equation*}
					and
		$$\norm{\psi_{\solg}}{\PSHM{m}(\Gp)}\leq C \norm{\solg}{\VSHM{m-3/2}_T(\Gamma)}$$	
		for $m=2,\ldots, \ell+2$, 
		where the constant $C>0$ depends only on $\Gp$, $\varepsilon, \zeta$ and $\ell$.
\end{lemma}
\begin{fatproof}
		We notice that $\diverg_{\Gamma}\zeta\nabla_{\Gamma}\psi_{\solg}=-i\diverg_{\Gamma}\solg$ is an elliptic problem on $\Gamma$, hence elliptic regularity theory asserts the existence of a function $\varphi\in\SHM{\ell+3/2}(\Gamma)$ which satisfies $\diverg_{\Gamma}\zeta\nabla_{\Gamma}\varphi=-i\diverg_{\Gamma}\solg$. The sought function $\psi_{\solg}$ can then be obtained as the solution of the Poisson problem
\begin{equation*}
				\begin{alignedat}{2}
					\diverg\varepsilon\nabla\psi_{\solg} &=0 \quad &&{\rm in} \quad \Omega, \\
				\psi_{\solg} & = \varphi \quad &&{\rm on} \quad \Gamma.
				\end{alignedat}
					\end{equation*}
					Due to \cite[Lem.~4.3]{MaxwellMyself} we have $\psi_{\solg}\in\PSHM{\ell+2}(\Gp)$ together with the estimate
					$\norm{\psi_{\solg}}{\PSHM{m}(\Gp)}\leq C \norm{\solg}{\VSHM{m-3/2}_T(\Gamma)}$ for all $m=2,\ldots,\ell+2$. Furthermore, by construction there holds $\diverg_{\Gamma}\zeta\nabla_{\Gamma}\psi_{\solg} = -i\diverg_{\Gamma}\solg$, hence the proof is complete.	
\end{fatproof}

Lemma~\ref{laxmilgram} asserts that the solution map of the auxiliary Maxwell problem $\solimpaux:\VSL(\Omega)\times\VSL_T(\Gamma)\rightarrow\HXI$ defined by $(\solf,\solgi)\mapsto\soluaux$ is well-defined and bounded. Furthermore, under additional assumptions on $\solf$ and $\solgi$, Lemma~\ref{regshift2} asserts additional regularity of $\solimpaux\left(\solf,\solgi\right)$. 
The following result combines this fact with the findings from Lemma~\ref{freqoplemma} and Proposition~\ref{freqopbdrprop}. We recall the operators $\Higeom$ and $\Higamma$ from Lemma~\ref{freqoplemma} and Proposition~\ref{freqopbdrprop} and remember that $\oldlambda > 1$ is a parameter that we will adjust to our needs later.

\begin{lemma}\label{contraction} 
		Let $\Gp = \geom$ be an $\ana$-partition and suppose that $\mu^{-1},\varepsilon\in\anatenspw{\Gp}$ satisfy \eqref{coercivemu}-\eqref{coerciveeps} and that $\zeta\in\anagamma$ satisfies \eqref{coercivezeta} and \eqref{zetaproperty}. Furthermore, let $k\in\Coone$ as well as $\solf\in\PVSHM{2}(\Gp)$ with $\diverg\solf = 0$ in $\Omega$ be given and suppose that $\solgi\in\VSHM{3/2}_T(\Gamma)$ and $\diverg_{\Gamma}\solgi = 0$. 

		Then, there holds
		\begin{align*}
				|k|^2\norm{\solimpaux\left(\Higeom\solf, \Higamma\solgi\right)}{\PVSHM{2}(\Gp)}&\leq C\oldlambda^{-1}\left(\norm{\solf}{\PVSHM{2}(\Gp)}+|k|\norm{\solgi}{\VSHM{3/2}_T(\Gamma)}\right)
		\end{align*}
		where the constant $C>0$ depends only on $\Gp$, $\mu^{-1}$, $\varepsilon$ and $\zeta$.
\end{lemma}

\begin{fatproof}
		According to Lemma~\ref{regshift2} we have
		\begin{align*}
				&|k|^2\norm{\solimpaux\left(\Higeom\solf, \Higamma\solgi\right)}{\PVSHM{2}(\Gp)} \\
				&\hskip 1cm \leq C|k|^2\left( 
				\norm{\Higeom\solf}{\VSL(\Omega)}+|k|^{1/2}\norm{\Higamma\solgi}{\VSL_T(\Gamma)}+\norm{\Higamma\solgi}{\VSHM{1/2}_T(\Gamma)}+|k|^{-1}\norm{\Higeom\solf\cdot\soln}{\SHM{1/2}(\Gamma)}\right).
		\end{align*}
		Due to Lemma~\ref{freqoplemma}, Proposition~\ref{freqopbdrprop} and the fact that $\diverg\solf=0$ implies $\diverg\Higeom\solf=0$ there holds
		\begin{align*}
				|k|^2\norm{\Higeom\solf}{\VSL(\Omega)}+|k|\norm{\Higeom\solf\cdot\soln}{\SHM{1/2}(\Gamma)} &\leq 
		|k|^2\norm{\Higeom\solf}{\VSL(\Omega)}+|k|\norm{\Higeom\solf}{\PVSHM{1}(\Gp)}  \\
																										   &\leq C\oldlambda^{-1}\norm{\solf}{\PVSHM{2}(\Gp)}
		\end{align*}
		as well as 
		\begin{align*}
				|k|^{5/2}\norm{\Higamma\solgi}{\VSL_T(\Gamma)}+|k|^2\norm{\Higamma\solgi}{\VSHM{1/2}_T(\Gamma)}\leq C\oldlambda^{-1}|k|\norm{\solgi}{\VSHM{3/2}_T(\Gamma)}.
		\end{align*}
		This concludes the proof. 
		
\end{fatproof}

In the subsequent proof of Theorem~\ref{Mainresult1} we use the abbreviations
\begin{align*}
		\curlcurl(\solv) &:= \curl\mu^{-1}\curl\solv-k^2\varepsilon\solv, \\
		\impedanz(\solv) &:= (\mu^{-1}\curl\solv)\times\soln -ik\zeta\solv_T.
\end{align*}
We recall the solution map $\solimp:\VSL(\Omega)\times\VSL_T(\Gamma)\rightarrow\HXI$ of the Maxwell problem \eqref{Maxwellorig}, which is well-defined for all $k\in\setu\subseteq\Coone$ provided $\setu\subseteq\Coone$ is such that Assumption~\ref{assumption1} holds. 
With these definitions and the preliminary auxiliary results we are finally able to present a proof of Theorem~\ref{Mainresult1}.

\bigskip 

\begin{fatproofmod}{Theorem~\ref{Mainresult1}} 
		We follow the lines of \cite[Thm.~7.3]{MaxwellImpedanceMelenk} and divide the proof into two steps. In the first step we reduce the data $\solf$ and $\solgi$ to divergence-free data $\solf'$ and $\solgi'$, and in the second step we present the main argument to conclude the proof. 

		\smallskip

		{\bf Step 1:} Our aim is to reduce the given data $\solf$ and $\solgi$ to divergence-free data $\solf'$ and $\solgi'$. To that end, let $\psi_{\solf}$ be the solution of
		\begin{alignat*}{2}
				\diverg\varepsilon\nabla\psi_{\solf} &= \diverg\solf \quad &&{\rm in}\ \Omega, \\	
				\psi_{\solf} &= 0 \quad &&{\rm on}\ \Gamma.
		\end{alignat*}
		According to \cite[Lem.~4.3]{MaxwellMyself} there holds $\psi_{\solf}\in\PVSHM{2}(\Gp)$ together with $\norm{\psi_{\solf}}{\PVSHM{2}(\Gp)}\leq C\norm{\diverg\solf}{\VSL(\Omega)}$.
		Then, with $\psi_{\solgi}$ being the function from Lemma~\ref{liftg} associated to $\solgi$, the vector field $\solu':=\solu+k^{-2}\nabla\psi_{\solf}+k^{-1}\nabla\psi_{\solgi}$ satisfies
		\begin{alignat*}{3}
				\curlcurl(\solu') & = \solf-\varepsilon\nabla\psi_{\solf}-k\varepsilon\nabla\psi_{\solgi} &&:= \solf' \quad && {\rm in}\ \Omega,\\
				\impedanz(\solu') &= \solgi-i\zeta\nabla_{\Gamma}\psi_{\solgi}&&:=\solgi'\quad && {\rm on}\ \Gamma, 
		\end{alignat*}
		and we observe $\diverg\solf' = 0$ as well as $\diverg_{\Gamma}\solgi' = 0$. Moreover, there hold the estimates
		\begin{align*}
				\norm{\solf'}{\PVSHM{1}(\Gp)}&\leq C\left(\norm{\solf}{\PVSHM{1}(\Gp)}+|k|\norm{\solgi}{\VSHM{1/2}_T(\Gamma)}\right), \\
				\norm{\solgi'}{\VSHM{1/2}_T(\Gamma)}&\leq C\norm{\solgi}{\VSHM{1/2}_T(\Gamma)},
		\end{align*}
		where the constant $C>0$ depends only on $\Gp$, $\varepsilon$ and $\zeta$.

		\medskip 

		{\bf Step 2:} The goal of the second step is to decompose $\solu'$ into $\solu' = \soluH+\soluA$ where $\soluH\in\PVSHM{2}(\Gp)$ and $\soluA$ is piecewise analytic. To that end we define the vector fields $\solf_0:=\solf'$ and $\solg_0:=\solgi'$, and for all $\ell\in\N_0$ 
		\begin{equation}\label{recursivedef}
		\begin{split}
				\solu_{\ell+1}&:=\solimpaux\left(\Higeom\solf_{\ell}, \Higamma\solg_{\ell}\right), \\
				\solw_{\ell+1} &:= \solimp\left(\Lowgeom\solf_{\ell}, \Logamma\solg_{\ell}\right),  \\
				\solf'_{\ell+1} &:= (1-\sigma_{\varepsilon})k^2\varepsilon\solu_{\ell+1}, \\
				\solg'_{\ell+1} &:= (1-\sigma_{\zeta})ik\zeta\Pi_T\solu_{\ell+1}, \\
				\solf_{\ell+1} &:= \solf'_{\ell+1}-k\varepsilon\nabla\psi_{\solg'_{\ell+1}}, \\
				\solg_{\ell+1} &:= \solg'_{\ell+1}-i\zeta\nabla_{\Gamma}\psi_{\solg'_{\ell+1}},
		\end{split}
\end{equation}
		where $\psi_{\solg'_{\ell+1}}$ is the function from Lemma~\ref{liftg} associated to $\solg'_{\ell+1}$. 
		Due to $\diverg\solf_0 = 0$ and $\diverg\varepsilon\nabla\psi_{\solg'_{\ell+1}}=0$ as well as $\diverg\Higeom\solf_{\ell} = 0$, we notice that $\diverg\solf_{\ell}=0$ for all $\ell\in\N_0$. Furthermore, by construction we have $\diverg_{\Gamma}\solg_{\ell}=0$ for all $\ell\in\N_0$. 

		With these preliminaries we employ the definition of $\solg_1$ and Lemma~\ref{liftg} to estimate 
		\begin{align*}
		\norm{\solf_1}{\PVSHM{2}{(\Gp)}}+|k|\norm{\solg_1}{\VSHM{3/2}_T(\Gamma)}\leq C|k|^2\norm{\solu_1}{\PVSHM{2}(\Gp)},
		\end{align*}
		and together with Lemma~\ref{regshift2} we get
		\begin{align*}
				|k|^2\norm{\solu_1}{\PVSHM{2}(\Gp)}&\leq C|k|^2\left(\norm{\Higeom\solf'}{\VSL(\Omega)}+|k|^{1/2}\norm{\Higamma\solgi'}{\VSL_T(\Gamma)}+\revision{\norm{\Higamma\solgi'}{\VSHM{1/2}_T(\Gamma)}}+|k|^{-1}\norm{\Higamma\solf'\cdot\soln}{\SHM{1/2}(\Gamma)}\right) \\
												   &\leq C|k|^2\left(|k|^{-1}\norm{\solf'}{\PVSHM{1}(\Gp)}+\norm{\solgi'}{\VSHM{1/2}_T(\Gamma)}\right) \\
												   &\leq C|k|\left(\norm{\solf}{\PVSHM{1}(\Gp)}+|k|\norm{\solgi}{\VSHM{1/2}_T(\Gamma)}\right).
		\end{align*}
		We conclude
		\begin{align}\label{f1impedance}
			\norm{\solf_1}{\PVSHM{2}{(\Gp)}}+|k|\norm{\solg_1}{\VSHM{3/2}_T(\Gamma)}
												   \leq C|k|\left(\norm{\solf}{\PVSHM{1}(\Gp)}+|k|\norm{\solgi}{\VSHM{1/2}_T(\Gamma)}\right),
		\end{align}
		which provides an estimate for $\solf_1$ and $\solg_1$. Next, we aim for similar estimates for $\solf_{\ell}$ and $\solg_{\ell}$ for $\ell>1$. We notice that due to $\solf_1\in\PVSHM{2}(\Gp)$ and $\solg_1\in\VSHM{3/2}_T(\Gamma)$, we may employ Lemma~\ref{contraction} to get better\footnote{\revision{By} better we mean better in terms of powers of $k$.} regularity estimates for $\solf_{\ell}$ and $\solg_{\ell}$ in the case $\ell>1$. 

	In fact, by definition of $\solf_{\ell+1}$,
		\begin{align*}
				\norm{\solf_{\ell+1}}{\PVSHM{2}(\Gp)}\leq C|k|^2\norm{\solu_{\ell+1}}{\PVSHM{2}(\Gp)}
		\end{align*}
		for all $\ell\geq 1$, and due to $\diverg_{\Gamma}\solg_{\ell+1} = 0$ and \cite[Lem.~3.6]{MaxwellMyself} we have for all $\ell\geq 1$
		\begin{align}\label{f2impedance}
				\begin{split}
						|k|\norm{\solg_{\ell+1}}{\VSHM{3/2}_T(\Gamma)}&\leq C|k|\norm{\curl_{\Gamma}\zeta^{-1}\solg_{\ell+1}}{\SHM{1/2}(\Gamma)}\\
																	  &\leq C|k|^2\norm{\curl\solu_{\ell+1}}{\PVSHM{1}(\Gp)}\\
																	  &\leq C|k|^2\norm{\solu_{\ell+1}}{\PVSHM{2}(\Gp)}.
		\end{split}
		\end{align}
		 Combining \eqref{f1impedance} and \eqref{f2impedance} and applying Lemma~\ref{contraction} yields 
		\begin{align}\label{f25impedance}
				\begin{split}
					\norm{\solf_{\ell+1}}{\PVSHM{2}(\Gp)}+|k|\norm{\solg_{\ell+1}}{\VSHM{3/2}_T(\Gamma)}&\leq C|k|^2\norm{\solu_{\ell+1}}{\PVSHM{2}(\Gp)}\\
																									&\leq C\oldlambda^{-1}\left(\norm{\solf_{\ell}}{\PVSHM{2}(\Gp)}+|k|\norm{\solg_{\ell}}{\VSHM{3/2}_T(\Gamma)}\right)
				\end{split}
		\end{align}
		for all $\ell\geq 1$ and a constant $C>0$ independent of $\ell$. By induction, \eqref{f25impedance} \revtwo{implies for $\ell\geq 1$}
		\begin{align}\label{f3impedance}
				\begin{split}
						\norm{\solf_{\ell+1}}{\PVSHM{2}(\Gp)}+|k|\norm{\solg_{\ell+1}}{\VSHM{3/2}_T(\Gamma)}&\leq 
				C|k|^2\norm{\solu_{\ell+1}}{\PVSHM{2}(\Gp)}\\
																											&\leq (\revtwo{C'}\oldlambda^{-1})^{\ell}\left(\norm{\solf_1}{\PVSHM{2}(\Gp)}+|k|\norm{\solg_1}{\VSHM{3/2}_T(\Gamma)}\right)
				\end{split}
					\end{align}
					\revtwo{for constant $C'>0$ that is independent of $k$ and $\ell$.} The above estimates allow us to control $\solf_{\ell}$ and $\solg_{\ell}$ by $\solf_1$ and $\solg_1$. Furthermore, by choosing $\oldlambda>\revtwo{C'}$ there holds $q:=\revtwo{C'}\oldlambda^{-1}<1$, hence \eqref{f3impedance} and \eqref{f1impedance} prove
		\begin{align}\label{f4impedance}
				\begin{split}
						\norm{\solf_{\ell+1}}{\PVSHM{2}(\Gp)}+|k|\norm{\solg_{\ell+1}}{\VSHM{3/2}_T(\Gamma)}&\leq C|k|^2\norm{\solu_{\ell+1}}{\PVSHM{2}(\Gp)}\\
																											&\leq C|k|q^{\ell}\left(\norm{\solf}{\PVSHM{1}(\Gp)}+|k|\norm{\solgi}{\VSHM{1/2}_T(\Gamma)}\right)
				\end{split}
		\end{align}
		for all $\ell\in\N_0$ for some $q\in (0,1)$ \revtwo{and a constant $C>0$ that is independent of $k$ and $\ell$}.

		\medskip

		While \eqref{f4impedance} provides an estimate for the growth of $\solf_{\ell}$ and $\solg_{\ell}$, we still lack an estimate for $\solw_{\ell}$. 
		According to Lemma~\ref{freqoplemma} and Proposition~\ref{freqopbdrprop} we have
		\begin{align*}
				\forall m\in\N_0:\ \norm{\Lowgeom\solf_{\ell}}{\PVSHM{m}(\Gp)}&\leq \omega M^{m}(m+|k|)^m\norm{\solf_{\ell}}{\VSL(\Omega)}, \\
				\Logamma\solg\in\anagammavecspecial{\Gamma}{\omega\norm{\solg_{\ell}}{\VSL_T{\Gamma}}}{M}, 
		\end{align*}
		where the constants $\omega>0$ is independent of $k$ and $\oldlambda$ and $M>0$ is independent of $k$, but may depend on $\oldlambda$.\footnote{The dependence of $\oldlambda$ is not a problem since $\oldlambda$ was already selected in \eqref{f4impedance}.}
		We appeal to \cite[Thm.~2.14]{MaxwellMyself2} and get
		\begin{align}\label{f5impedance}
				\forall m\in\N_0:\ \norm{\solw_{\ell+1}}{\PVSHM{m}(\Gp)}&\leq CC_k L^m(m+|k|)^m,
		\end{align}
	 where
		$C>0$ depends only on $\Gp$ and $\omega$, the constant $L>0$ depends only on $\Gp$, $\mu^{-1}$, $\varepsilon$, $\zeta$, $M$ and the previously fixed $\oldlambda$, and where
		\begin{align*}
				C_k :=\norm{\solw_{\ell+1}}{\VSL(\Omega)}+|k|^{-1}\norm{\curl\solw_{\ell+1}}{\VSL(\Omega)}+|k|^{-2}\norm{\solf_{\ell}}{\VSL(\Omega)}+|k|^{-3/2}\norm{\solg_{\ell}}{\VSL_T(\Gamma)}. 
		\end{align*}
		Moreover, by \eqref{rhodef} and \eqref{f4impedance} we may estimate
		\begin{align*}
				\norm{\solw_{\ell+1}}{\VSL(\Omega)}+|k|^{-1}\norm{\curl\solw_{\ell+1}}{\VSL(\Omega)} &\leq |k|^{-1}\rho(k)\left(\norm{\Lowgeom\solf_{\ell}}{\VSL(\Omega)}+\norm{\Logamma\solg_{\ell}}{\VSL_T(\Gamma)}\right)\\
																									 &\leq C|k|^{-1}\rho(k)\left(\norm{\solf_{\ell}}{\VSL(\Omega)}+\norm{\solg_{\ell}}{\VSL_T(\Gamma)}\right)\\
																									 &\leq C\rho(k)q^{\ell}\left(\norm{\solf}{\PVSHM{1}(\Gp)}+|k|\norm{\solgi}{\VSHM{1/2}_T(\Gamma)}\right),
		\end{align*}
		and
		\begin{align*}
				|k|^{-2}\norm{\solf_{\ell}}{\VSL(\Omega)}+|k|^{-3/2}\norm{\solg_{\ell}}{\VSL_T(\Gamma)}\leq Cq^{\ell}(|k|^{-1}+|k|^{-3/2})\left(\norm{\solf}{\PVSHM{1}(\Gp)}+|k|\norm{\solgi}{\VSHM{1/2}_T(\Gamma)}\right).
		\end{align*}
		Combining these estimates with the definition of $C_k$ and \eqref{f5impedance} yields
		\begin{align}\label{f6impedance}
				\forall m\in\N_0:\ \norm{\solw_{\ell+1}}{\PVSHM{m}(\Gp)}&\leq Cq^{\ell}\left[\rho(k)+|k|^{-1}\right] \left(\norm{\solf}{\PVSHM{1}(\Gp)}+|k|\norm{\solgi}{\VSHM{1/2}_T(\Gamma)}\right) L^m(m+|k|)^m
		\end{align}
		for some $q\in (0,1)$. This provides the sought estimate for $\solw_{\ell+1}$.
		
	\bigskip	

	We now enter the last stage of the proof. Recalling \eqref{recursivedef}, we define the vector fields $\soluH'$ and $\soluA$ via
		\begin{align*}
				\soluH':=\sum_{\ell=1}^{\infty}\solu_{\ell}, \quad \soluA := \sum_{\ell=1}^{\infty}\solw_{\ell}, \quad \nabla\psi:=\sum_{\ell=1}^{\infty}\nabla\psi_{\solg'_{\ell}}, \quad \soluH:=\soluH'-k^{-1}\nabla\psi.
		\end{align*}
		From Lemma~\ref{liftg}, the definition of $\solg'_{\ell}$ in \eqref{recursivedef} and a trace inequality we obtain for all $\ell\geq 1$
		\begin{align*}
				\norm{\nabla\psi_{\solg'_{\ell}}}{\PVSHM{2}(\Gp)}&\leq C\norm{\solg'_{\ell}}{\VSHM{3/2}_T(\Gamma)} \\
																 &\leq C|k|\norm{\solu_{\ell}}{\PVSHM{2}(\Gp)}.
		\end{align*}
		Hence, \eqref{f4impedance} and a geometric series argument prove
		\begin{align*}
				\norm{\nabla\psi}{\PVSHM{2}(\Gp)}&\leq C|k|\sum_{\ell=1}^{\infty}\norm{\solu_{\ell}}{\PVSHM{2}(\Gp)}\\
												 &\leq C\left(\norm{\solf}{\PVSHM{1}(\Gp)}+|k|\norm{\solgi}{\VSHM{1/2}_T(\Gamma)}\right),
		\end{align*}
		and therefore
		\begin{align}\label{f7impedance}
				\norm{\soluH}{\PVSHM{2}(\Gp)}\leq C|k|^{-1}\left(\norm{\solf}{\PVSHM{1}(\Gp)}+|k|\norm{\solgi}{\VSHM{1/2}_T(\Gamma)}\right).
		\end{align}
		Similarly, \eqref{f6impedance} and a geometric series argument show
		\begin{align}\label{f8impedance}
				\forall m\in\N_0:\ \norm{\soluA}{\PVSHM{m}(\Gp)}\leq C\left[\rho(k)+|k|^{-1}\right] \left(\norm{\solf}{\PVSHM{1}(\Gp)}+|k|\norm{\solgi}{\VSHM{1/2}_T(\Gamma)}\right)L^m(m+|k|)^m.
		\end{align}

		We notice that \eqref{f7impedance}-\eqref{f8impedance} are the claimed estimates for $\soluH$ and $\soluA$. It remains to show that indeed $\solu'=\soluH+\soluA$, where $\solu':=\solu+k^{-2}\nabla\psi_{\solf}+k^{-1}\nabla\psi_{\solgi}$ is the modified solution constructed in Step 1. To that end we use the unique solvability assumption from Assumption~\ref{assumption1} which implies that $\solu'=\soluH+\soluA$ is equivalent to $\curlcurl\left(\soluH+\soluA\right) = \curlcurl(\solu')$ and $\impedanz\left(\soluH+\soluA\right)=\impedanz(\solu')$. 

		Indeed, a calculation shows
	\begin{align*}
			\curlcurl\left(\soluH+\soluA\right) = \sum_{\ell=1}^{\infty}\curlcurl\left(\solu_{\ell}-k^{-1}\nabla\psi_{\solg'_{\ell}}+\solw_{\ell}\right) = \sum_{\ell=1}^{\infty}\left[\solf_{\ell-1}-\solf_{\ell}\right]= \solf',
	\end{align*}
	as well as 
	\begin{align*}
			\impedanz\left(\soluH+\soluA\right) = \sum_{\ell=1}^{\infty}\impedanz\left(\solu_{\ell}-k^{-1}\nabla\psi_{\solg'_{\ell}}+\solw_{\ell}\right) = \sum_{\ell=1}^{\infty}\left[\solg_{\ell-1}-\solg_{\ell}\right] = \solg',	
	\end{align*}
	hence $\solu' = \soluH+\soluA$ by Assumption~\ref{assumption1}. This concludes the proof.
	
\end{fatproofmod}

	\section{Discretization}\label{discretizationsec}

After having discussed the regularity properties of Maxwell's equations we turn our attention to the finite-element approximations of a solution $\solu$. The aim of this intermediate section is to recall the definition of conforming finite-element spaces based on \ned-type I elements \cite{nedelectype1} along with their properties. Throughout this section we closely follow \cite[Sec. 8]{MaxwellImpedanceMelenk}.

\subsection{Regular triangulations and \ned-type I elements}

Let $\Gp=\geom$ be an $\ana$-partition in the sense of Definition~\ref{partitiondef} and let $\Tref$ be the (open) reference tetrahedron given by the interior of the convex hull of the origin and \revtwo{three distinct points}. A set $\T$ consisting of finitely many open subsets of $\R^3$ (which we call elements) is called a regular and shape regular triangulation of $\Gp$ if 
\begin{itemize}
		\item The elements $T\in\T$ satisfy $\overline{\Omega} = \bigcup_{T\in\T}\overline{T}$,
		\item For all elements $T\in\T$ and all (by definition open) subdomains $\Gp_i$ there holds either $T\subseteq\Gp_i$ or $T\subseteq\Omega\setminus\overline{\Gp_i}$,
		\item For every $T\in\T$ there exists a $\CM{1}$-diffeomorphism $F_T:\overline{\Tref}\rightarrow\overline{T}$.
		\item For two distinct elements $T,T'\in\T$, the intersection $\overline{T}\cap\overline{T'}$ is either empty, a common vertex, a common edge or a common face, where vertices, edges and faces of $T$ and $T'$ are defined as the images of the corresponding quantities of $\Tref$ under $F_T$ and $F_{T'}$, respectively. 

				In addition, we suppose that for elements $T,T'$ sharing an edge or a face, the respective maps $F_T$ and $F_{T'}$ are compatible, meaning that if $T,T'$ share an edge (i.e., $F_T(e) = F_{T'}(e')$ for edges $e,e'$ of $\Tref$) or share a face (i.e., $F_T(f)=F_{T'}(f')$ for faces $f,f'$ of $\Tref$), then $F_T\circ F_{T'}^{-1}:e'\rightarrow e$ or $F_T\circ F_{T'}^{-1}:f'\rightarrow f$ are affine bijections.
		\item For all $T\in\T$ there holds, with some shape-regularity constant $\gamma>0$,
				\begin{align*}
						h_T^{-1}\norm{F'_T}{\Linfty(\Tref)}+h_T\norm{(F'_T)^{-1}}{\Linfty(\Tref)}\leq \gamma,	
				\end{align*}
				where $h_T:=\operatorname{diam}(T)$.
\end{itemize}

We notice that the condition $T\subseteq\Gp_i$ or $T\subseteq\Omega\setminus\overline{\Gp_i}$ ensures that a regular and shape regular triangulation $\T$ resolves the subdomain interfaces $\interf_j$, meaning that every element belongs to exactly one subdomain. 

Intuitively, we may think of $\T$ as a finite union of tetrahedra, however due to the boundary $\Gamma$ and the subdomain interfaces $\interf_j$ being analytic surfaces, it is clear that at least the elements bordering $\Gamma$ and/or an interface $\interf_j$ must have curved edges and/or curved faces. 

However, in order to do a finite element analysis based on regular and shape regular triangulations $\T$ we need at least some control over the shapes and curvatures of the elements $T\in\T$. To that end we define the maximal mesh width $h:=\max\left\{h_T\ |\ T\in\T\right\}$ and suppose that all element maps $F_T$ can be written as a composition of an $h$-dependent affine scaling and an $h$-independent diffeomorphism. More precisely, we follow e.g. \cite[Sec.~8]{MaxwellImpedanceMelenk} and make the following assumption:

\begin{assumption}\label{Ttrafo}
$\T$ is a regular and shape regular triangulation of $\Gp$ such that every element map $F_T$ can be written as $F_T = S_T\circ A_T$, where $A_T$ is an affine diffeomorphism, $S_T$ is an analytic map that has an analytic inverse, and there exist parameters $C,\gamma>0$ independent of $T$ such that
\begin{align*}
		\forall T\in\T:\ & h_T^{-1}\norm{A'_T}{\Linfty(\Tref)}+h_T\norm{(A'_T)^{-1}}{\Linfty(\Tref)}+\norm{(S'_T)^{-1}}{\Linfty(\widetilde{T})}\leq C, \\
		\forall T\in\T:\ \forall \ell\in\N_0:\ &\norm{\nabla^{\ell}S_T}{\Linfty(\widetilde{T})}\leq C\gamma^{\ell}\ell^{\ell},
\end{align*}
where $h_T>0$ is the element diameter of $T$ and $\widetilde{T} := A_T(\Tref)$ is the image of $\Tref$ under $A_T$.
\end{assumption}

After having discussed the notion of regular and shape regular triangulations we may focus on finite element spaces based on \ned-type I elements of polynomial order $p\geq 0$: On the reference tetrahedron $\Tref$ we define
\begin{align*}
		\poly_p(\Tref)&:=\operatorname{span}\left\{{\bf x}^{\alpha}\ |\ |\alpha|\leq p\right\}, \\
		\nedelec(\Tref)& :=\left\{{\bf p}({\bf x})+{\bf x}\times{\bf q}({\bf x})\ |\ {\bf p},{\bf q}\in\left(\mathcal{P}_p(\Tref)\right)^3\right\}.
\end{align*}
For a regular and shape regular triangulation $\T$, the conforming spaces $S_{p+1}$ and $\nedelec(\Tref)$ are then defined by
\begin{align*}
	S_{p+1}(\T)&:=\left\{u\in\SHM{1}(\Omega)\ |\ u\vert_{T}\circ F_T\in\mathcal{P}_{p+1}(\Tref)\right\}, \\
	\nedelec(\T)&:=\left\{\solu\in\Hcurl\ |\ (F'_T)^T\solu\vert_T\circ F_T\in\nedelec(\Tref)\right\}.
\end{align*}

\begin{remark}\label{traforemark}
		For two bounded Lipschitz domains $\Omega_1$ and $\Omega_2$ and a diffeomorphism $\varphi:\Omega_2\rightarrow\Omega_1$ with $F:=\varphi'$ and $J:=\operatorname{det}F$, let $\wu:=F^T(\solu\circ\varphi)$ and $\solucheck:=JF^{-1}(\solu\circ\varphi)$ denote the covariant and contravariant transformation of a vector field $\solu$ on $\Omega_2$, respectively.

		With this notation, we observe that $\nedelec(\T)$ consists of all vector fields $\solu\in\Hcurl$ whose element-wise covariant transformations belong to $\nedelec(\Tref)$.
\end{remark}

We notice that $\nedelec(\T)\subseteq\HXI$, hence we may define conforming spaces of discrete vector fields and discrete scalar potentials by
\begin{align*}
		\Xh:=\nedelec(\T) \quad {\rm and}\quad \Sh:=S_{p+1}(\T).
\end{align*}
There holds the exact sequence property \cite[Lem.~5.28]{BookMonk}
\begin{align}\label{exactseq}
		\Sh\xrightarrow{\ \ \nabla\ \ }\Xh \xrightarrow{\ \ \curl\ \ }\curl\Xh.
\end{align}

The discrete variational formulation of \eqref{Maxwellorig} then reads as: Find $\soluh\in\Xh$ such that
\begin{align}\label{discreteproblem}
		\begin{split}
				\forall\solvh\in\Xh:\ \SCP{\mu^{-1}\curl \soluh}{\curl\solvh}{\VSL(\Omega)}-k^2\SCP{\varepsilon\soluh}{\solvh}{\VSL(\Omega)}-i&k\SCP{\zeta\Pi_T\soluh}{\Pi_T\solvh}{\VSL_T(\Gamma)} \\
		&= \SCP{\solf}{\solvh}{\VSL(\Omega)}+\SCP{\solgi}{\Pi_T\solvh}{\VSL_T(\Gamma)}.
\end{split}
\end{align}

\subsection{Element-by-element approximation operators}

Following \cite{MaxwellImpedanceMelenk} we rely on special $hp$-approximation operators that were constructed in \cite{MaxwellTransparentMelenk, MaxwellImpedanceMelenk, OperatorsRojik}. For every element $T\in\T$ we define
\begin{align*}
		\HcurlT{1}&:=\left\{\solu\in\VSHM{1}(T)\ |\ \curl\solu\in\VSHM{1}(T)\right\}, \\
		\HdivT{1}&:=\left\{\solu\in\VSHM{1}(T)\ |\ \diverg\solu\in\SHM{1}(T)\right\}
\end{align*}
and set 
\begin{align*}
		\Hcurlth{1}&:=\Hcurl\cap\prod_{T\in\T}\HcurlT{1}, \\
		\Hdivth{1}&:=\Hdiv\cap\prod_{T\in\T}\HdivT{1}.
\end{align*}

Henceforth, we call an operator $\Pi:\SHM{m}(\Omega)\rightarrow S_{p+1}(\T)$ an element-by-element operator if there exists an operator $\Piref:\SHM{m}(\Tref)\rightarrow\poly_p(\Tref)$ such that $(\Pi u)\vert_{T}\circ F_T = \Piref(u\circ F_T)$ for all $u\in\SHM{m}(\Omega)$ and all elements $T\in\T$.

The following result \cite[Cor.~8.4]{MaxwellImpedanceMelenk} asserts the existence of element-by-element operators with optimal $hp$-approximation properties.

\begin{proposition}\label{defprigradhp}
		Let $m\geq 2$. For every $p\geq m-2$ there exists an element-by-element operator $\Pigradhp$ which has the mapping property $\Pigradhp:\SHM{m}(\Omega)\rightarrow S_{p+1}(\T)$ and features the approximation properties
		\begin{align*}
				\forall T\in\T:\ \norm{\varphi-\Pigradhp\varphi}{\SL(T)}+\frac{h_T}{p+1}\norm{\varphi-\Pigradhp\varphi}{\SHM{1}(T)}&\leq C'\left(\frac{h_T}{p+1}\right)^m\norm{\varphi}{\SHM{m}(T)}, \\
				\forall T\in\T:\ \norm{\varphi-\Pigradhp\varphi}{\SHM{1}(\partial T)}&\leq C'\left(\frac{h_T}{p+1}\right)^{m-3/2}\norm{\varphi}{\SHM{m}(T)}
		\end{align*}
		for a constant $C'>0$ that depends only on $m$ and the parameters $C,\gamma>0$ from Assumption~\ref{Ttrafo}.
\end{proposition}

Similarly to the scalar case, we call an operator $\Pi:\VSHM{m}(\Omega)\rightarrow\Xh$ an element-by-element operator, if there exists an operator $\Piref:\VSHM{m}(\Tref)\rightarrow\nedelec(\Tref)$ such that $\widehat{\Pi\solu\vert_T} = \Piref(\widehat{\solu\vert_T})$ for all $\solu\in\VSHM{m}(\Omega)$ and all $T\in\T$, where $\widehat{\Pi\solu\vert_T}$ and $\widehat{\solu\vert_T}$ are the covariant transformations of $\Pi\solu\vert_T$ and $\solu\vert_T$, respectively (cf. Remark~\ref{traforemark}). \revision{The following result \cite[Lem.~8.5]{MaxwellImpedanceMelenk} generalizes the previous proposition to the case of vector fields.} We recall that $h:=\max_{T\in\T} h_T$ denotes the maximal mesh width.

\begin{proposition}\label{defpicurlhp}
		Let $m\geq 2$. For every $p\geq m-1$ there exists an element-by-element operator $\Picurlhp$ which has the mapping property $\Picurlhp:\VSHM{m}(\Omega)\rightarrow\Xh$ and features the approximation properties
\begin{align*}
		\forall T\in\T:\ \norm{\solu-\Picurlhp\solu}{\VSL(T)}+\frac{h_T}{p+1}\norm{\solu-\Picurlhp\solu}{\VSHM{1}(T)}&\leq C'\left(\frac{h_T}{p+1}\right)^{m-1}\norm{\solu}{\VSHM{m}(T)}, \\
				\forall T\in\T:\ \norm{\solu-\Picurlhp\solu}{\VSL(\partial T)}&\leq C'\left(\frac{h_T}{p+1}\right)^{m-1/2}\norm{\solu}{\VSHM{m}(T)}
		\end{align*}
		for a constant $C'>0$ that depends only on $m$ and the parameters $C,\gamma>0$ from Assumption~\ref{Ttrafo}.

		Moreover, let $C_1,M>0$ be arbitrary and suppose that $\solu\in\anavecCM{\Omega}{\omega_\solu}{M}$ for some $\omega_{\solu}>0$. Then, there exist constants $C_2, \sigma>0$ depending only on $C_1, M$ and the parameters $C,\gamma>0$ from Assumption~\ref{Ttrafo} such that
		\begin{align*}
		h+|k|h/p\leq C_1
		\end{align*}
		implies
		\begin{align*}
		\norm{\solu-\Picurlhp\solu}{\HXIK}\leq C_2 \omega_{\solu}|k|\left(\left(\frac{h}{h+\sigma}\right)^p+\left(\frac{|k|h}{\sigma p}\right)^p\right).
		\end{align*}
\end{proposition}

Although the operator $\Picurlhp$ from Proposition~\ref{defpicurlhp} is nice in the sense that it is an element-by-element operator with optimal $hp$-approximation properties, it has the drawback that is not a projection onto $\Xh$. The subsequent result provides an operator $\Picurlcom$ which is a projection operator onto $\Xh$ at the expense of not having $hp$-optimal approximation properties in $\VSHM{1}$. 

\smallskip

Similarly as before we call an operator $\Pi:\Hcurlth{1}\rightarrow\Xh$ or $\Pi:\Hdivth{1}\rightarrow\curl\Xh$ an element-by-element operator if there exists an operator $\Piref:\VSHM{1}(\Tref)\rightarrow\nedelec(\Tref)$ or $\Piref:\VSHM{1}(\Tref)\rightarrow\curl\nedelec(\Tref)$ such that $\widehat{\Pi\solu\vert_T} = \Piref(\widehat{\solu\vert_T})$ or $\widecheck{\Pi\solu\vert_T} = \Piref(\widecheck{\solu\vert_T})$ for all $T\in\T$ and all $\solu\in\Hcurlth{1}$ or all $\solu\in\Hdivth{1}$, respectively.

With this notation there holds the following result \cite[Lem.~8.8]{MaxwellImpedanceMelenk}:

\begin{proposition}\label{Ipprop}
For every $p\in\N_0$ there exist element-by-element operators $\Picurlcom:\Hcurlth{1}\rightarrow\Xh$ and $\Pidivcom:\Hdivth{1}\rightarrow\curl\Xh$ with the properties
\begin{itemize}
		\item $\Picurlcom\solvh = \solvh$ for all $\solvh\in\Xh$, that is, $\Picurlcom$ is a projection onto $\Xh$,
		\item For all $\solv\in\Hcurlth{1}$ there holds $\curl\Picurlcom\solv = \Pidivcom\curl\solv$, that is, $\Picurlcom$ and $\Pidivcom$ commute with respect to the curl operator.
\end{itemize}

\end{proposition}

We notice that $\Picurlcom$ is nice in the sense that it projects onto the discrete space $\Xh$, however, it does not feature the optimal $hp$-approximation properties of $\Picurlhp$. A part of the proof of Theorem~\ref{Mainresult2} below will be to compose the operators $\Picurlhp$ and $\Picurlcom$ in a clever way so that we can exploit both the $hp$-approximation properties of $\Picurlhp$ and the projection properties of $\Picurlcom$.

	\section{Wavenumber-explicit finite element analysis}\label{femanalysis}

The aim of this section is to derive a wavenumber-explicit convergence analysis of $hp$-FEM for Maxwell's equations in anisotropic and piecewise smooth media. In particular, we will show that under certain scale-resolution conditions the finite element method produces quasi-optimal discrete solutions. To this end we rely on the so-called Schatz argument \cite{Schatz} to estimate the finite-element error by the best-approximation error before employing the arguments from \cite[Sec.~9]{MaxwellImpedanceMelenk} to bound  the best-approximation error explicitly in $h$, $p$ and $k$.

\subsection{The Schatz argument}\label{Schatzsection}

For now, our goal is to show that if the discrete space $\Xh$ is ``sufficiently rich", the discrete approximation $\soluh$ of Maxwell's equations \eqref{Maxwellorig} is quasi-optimal in the sense that the error between $\solu$ and $\soluh$ can be controlled in terms of the error between $\solu$ and the best approximation of $\solu$ in $\Xh$. 

We suppose that we are interested in Maxwell's equations \eqref{Maxwellorig} at wavenumber $k\in\setu$ for some region $\setu\subseteq\Coone$. Furthermore, we suppose that $\mu^{-1}$ and $\zeta$ satisfy \eqref{coercivemu} and \eqref{zetaproperty} and that $\varepsilon$ and $\zeta$ satisfy the uniform coercivity conditions \eqref{ezcoerceps}-\eqref{ezcoerczeta} for all $k\in\setu$. We recall the sesquilinear form 
\begin{align*}
\bfa(\solu,\solv) = \SCP{\mu^{-1}\curl \solu}{\curl\solv}{\VSL(\Omega)}-k^2\SCP{\varepsilon\solu}{\solv}{\VSL(\Omega)}-ik\SCP{\zeta\solu_T}{\solv_T}{\VSL(\Gamma)}
\end{align*}
from \eqref{Adef}
and for $\solu,\solv\in\HXI$ we abbreviate
\begin{align*}
		\dualterm{\solu}{\solv}:=k^2\SCP{\varepsilon\solu}{\solv}{\VSL(\Omega)}+ik\SCP{\zeta\solu_T}{\solv_T}{\VSL_T(\Gamma)}.
\end{align*}
Due to \eqref{coercivemu} and \eqref{ezcoerceps}-\eqref{ezcoerczeta} we may estimate
\begin{align}\label{normestimate}
\forall\solu\in\HXI:\ \norm{\solu}{\HXIK}^2\leq \cco\left|A_k(\solu,\solu)+(1+\alpha_k)\dualterm{\solu}{\solu}\right|,
\end{align}
where $\alpha_k\in\Co$ is unimodular and $\cco>0$, see \eqref{weakcoercivity}.
Based on the sesquilinear form $\dualterm{\cdot}{\cdot}$, we define the quantity
\begin{align*}
		\delta_k(\solu):=\sup_{\solw_N\in\Xh\setminus\{{\bf 0}\}}\left(2\frac{|\dualterm{\solu}{\solw_N}|}{\norm{\solu}{\HXIK}\norm{\solw_N}{\HXIK}}\right)
\end{align*}
for all non-zero $\solu\in\HXI$, and $\delta_k({\bf 0}):=0$. With these definitions there holds the following lemma, which is an extension of \cite[Prop.~9.1]{MaxwellImpedanceMelenk}:
\begin{proposition}\label{quasioptimalprop}
		Suppose that $\mu$ and $\zeta$ satisfy \eqref{coercivemu} and \eqref{zetaproperty}, respectively, and that $\varepsilon$ and $\zeta$ satisfy the \revtwo{uniform} coercivity properties \eqref{ezcoerceps}-\eqref{ezcoerczeta} for all $k\in\setu\subseteq\Coone$.
		Under these assumptions, let $\solu\in\HXI$ and $\solu_N\in\Xh$ be such that the Galerkin orthogonality $$\forall\solv_N\in\Xh:\ \bfa(\solu-\solu_N,\solv_N)=0$$ holds true. Then, for all $k\in\setu$ there holds that $\delta_k(\solu-\solu_N)<\cco^{-1}$ implies
\begin{align*}
		\norm{\solu-\solu_N}{\HXIK}\leq C\cco\frac{1+\delta_k(\solu-\solu_N)}{1-\cco\delta_k(\solu-\solu_N)}\ \inf_{\solw_N\in\Xh}\norm{\solu-\solw_N}{\HXIK},
\end{align*}
where $\cco>0$ is the ellipticity constant from \eqref{coercivemu} and \eqref{ezcoerceps}-\eqref{ezcoerczeta}, and $C>0$ depends only on $\mu^{-1}, \varepsilon$ and $\zeta$.

That is, if $\delta_k(\solu-\solu_N)$ is sufficiently small, the discrete vector field $\solu_N$ is quasi-optimal in $\Xh$, and the quasi-optimality constant depends explicitly on $\delta_k(\solu-\solu_N)$.
\end{proposition}
\begin{fatproof}
		We set $\sole_N:=\solu-\solu_N$ and choose an arbitrary $\solw_N\in\Xh$. Due to \eqref{normestimate} and $\bfa(\sole_N,\solw_N-\solu_N)=0$ there holds 
\begin{align*}
		\norm{\sole_N}{\HXIK}^2&\leq \cco\left|A_k(\sole_N,\solu-\solw_N)+(1+\alpha_k)\dualterm{\sole_N}{\solu-\solw_N}+(1+\alpha_k)\dualterm{\sole_N}{\solw_N-\solu_N}\right| \\
								  &\leq C\cco\norm{\sole_N}{\HXIK}\norm{\solu-\solw_N}{\HXIK}+\cco\delta_k(\sole_N)\norm{\sole_N}{\HXIK}\norm{\solw_N-\solu_N}{\HXIK} \\
								  &\leq C\cco\norm{\sole_N}{\HXI}\norm{\solu-\solw_N}{\HXI}+\cco\delta_k(\sole_N)\norm{\sole_N}{\HXIK}\norm{\solu-\solw_N}{\HXIK} \\
								  &\hskip 9cm +\cco\delta_k(\sole_N)\norm{\sole_N}{\HXIK}^2,
\end{align*}
where we used the triangle inequality for the last step. Rearranging proves the claim.
\end{fatproof}

\subsection{Continuous and discrete Helmholtz decompositions}

The preceding lemma shows that the finite element solution $\solu_h$ is quasi-optimal if the error $\sole_h:=\solu-\solu_h$ is such that $\delta_k(\sole_h)$ is small. The final goal of this work is to estimate $\delta_k(\sole_h)$ solely in terms of the maximal mesh width $h$ of the underlying triangulation $\T$, the local polynomial degree $p$ of the discrete space $\Xh$ and the wavenumber $k$. In order to analyze $\delta_k(\sole_h)$ we follow \cite[Sec.~9]{MaxwellImpedanceMelenk} and start by studying discrete Helmholtz decompositions on $\Xh$. 

We suppose that $\mu^{-1},\varepsilon\in\anatenspw{\Gp}$ and $\zeta\in\anatens{\Gamma}$ as well as $\setu\subseteq\Coone$ satisfy Assumption~\ref{assumption2}, and for $k\in\setu$ and $\solu,\solv\in\HXI$ we consider the sesquilinear form
\begin{align*}
		\dualtermdual{\solu}{\solv}^H:=\overline{k}^2\SCP{\varepsilon^{H}\solu}{\solv}{\VSL(\Omega)}-i\overline{k}\SCP{\zeta^H\solu_T}{\solv_T}{\VSL_T(\Gamma)},
\end{align*}
where $\varepsilon^H$ and $\zeta^H$ denote the hermitian transposes of $\varepsilon$ and $\zeta$, respectively.

We observe that
\begin{align}\label{dualtermdual}
	\forall\solu,\solv\in\HXI:\	\dualterm{\solu}{\solv} = \overline{\dualtermdual{\solv}{\solu}^H},
\end{align}
and following \cite[Sec.~9]{MaxwellImpedanceMelenk} we define continuous and discrete wavenumber-dependent Helmholtz operators $\Pigradkdual:\HXI\rightarrow\nabla\Honeimp$ and $\Pigradkhdual:\Xh\rightarrow\nabla\Sh$ by
\begin{alignat*}{2}
		\dualtermdual{\Pigradkdual\solv}{\nabla\psi}^H &= \dualtermdual{\solv}{\nabla\psi}^H \qquad &&\revtwo{\forall\psi\in\Honeimp,} \\
		\dualtermdual{\Pigradkhdual\solv_N}{\nabla\psi_N}^H &= \dualtermdual{\solv_N}{\nabla\psi_N}^H \qquad &&\revtwo{\forall\psi_N\in\Sh}.
\end{alignat*}
Due to \eqref{ezcoerceps}-\eqref{ezcoerczeta} and the Lax-Milgram lemma these operators are well-defined for all $k\in\setu$, hence they give rise to stable continuous and discrete Helmholtz decompositions
\begin{equation}\label{discretehelmholtz}
		\begin{alignedat}{2}
				\solv &= \Pigradkdual\solv+\Picurlkdual\solv, \qquad&&\revtwo{\forall\solv\in\HXI},\\
				\solv_N &= \Pigradkhdual\solv_N+\Picurlkhdual\solv_N, \qquad&&\revtwo{\forall\solv_N\in\Xh}.
		\end{alignedat}
\end{equation}
From \eqref{dualtermdual} we observe that the operators $\Pigradkdual$ and $\Pigradkhdual$ satisfy
\begin{equation}\label{defprop}
		\begin{alignedat}{2}
				\dualterm{\nabla\psi}{\Pigradkdual\solv} &= \dualterm{\nabla\psi}{\solv} \qquad &&\revtwo{\forall\psi\in\Honeimp,} \\
				\dualterm{\nabla\psi_N}{\Pigradkhdual\solv_N} &= \dualterm{\nabla\psi_N}{\solv_N}\qquad &&\revtwo{\forall\psi_N\in\Sh}.
\end{alignedat}
\end{equation}
Hence, with $\sole_N:=\solu-\solu_N$ we may write for $\solw_N\in\Xh$
\begin{align*}
		\dualterm{\sole_N}{\solw_N} = \dualterm{\sole_N}{(\Picurlkhdual-\Picurlkdual)\solw_N}+\dualterm{\sole_N}{\Picurlkdual\solw_N}+\dualterm{\sole_N}{\Pigradkhdual\solw_N}
\end{align*}
and we note that due to Galerkin orthogonality there holds $\dualterm{\sole_N}{\Pigradkhdual\solw_N}=0$, that is
\begin{align}\label{decompositionconsistency}
		\dualterm{\sole_N}{\solw_N} = \underbrace{\dualterm{\sole_N}{(\Picurlkhdual-\Picurlkdual)\solw_N}}_{:=T_1}+\underbrace{\dualterm{\sole_N}{\Picurlkdual\solw_N}}_{:=T_2}.
\end{align}

\subsection{The commuting diagram argument}
The aim of this subsection is to provide an estimate of the term $T_1$ in \eqref{decompositionconsistency}. To that end, we define the norm
\begin{align*}
		\norm{\solu}{\knormplus}:=|k|\norm{\solu}{\VSL(\Omega)}+|k|^{1/2}\norm{\solu_T}{\VSL_T(\Gamma)}
\end{align*}
and observe
\begin{align}\label{codi1}
		\left| \dualterm{\sole_N}{(\Picurlkhdual-\Picurlkdual)\solw_N}\right|\leq C\norm{\sole_N}{\knormplus}\norm{(\Picurlkhdual-\Picurlkdual)\solw_N}{\knormplus}
\end{align}
for a constant $C>0$ depending only on $\varepsilon$ and $\zeta$. We proceed as in \cite[Sec.~9]{MaxwellImpedanceMelenk}: Due to \eqref{discretehelmholtz} we have $\curl\Picurlkdual\solw_N = \curl\Picurlkhdual\solw_N = \curl\solw_N$ for all $\solw_N\in\Xh$. Hence, with $\Picurlcom$ and $\Pidivcom$ being the \revtwo{projectors} from Proposition~\ref{Ipprop} there follows
\begin{align*}
		\curl\left(\Picurlkhdual\solw_N-\Picurlcom\Picurlkdual\solw_N\right) &= \curl\left(\Picurlkhdual\solw_N\right)-\Pidivcom\curl\left(\Picurlkdual\solw_N\right)\\
																			 &=\curl\solw_N-\curl\Picurlcom\solw_N \\
																			 &= \curl\solw_N-\curl\solw_N = 0.
\end{align*}
Due to the exact sequence property \eqref{exactseq} this shows $\Picurlkhdual\solw_N-\Picurlcom\Picurlkdual\solw_N = \nabla\xi_N$ for some $\xi_N\in\Sh$, and rewriting yields 
\begin{align*}
\left(\Picurlkhdual-\Picurlkdual\right)\solw_N = \nabla\xi_N+\left(\left(\Picurlcom-\operatorname{I}\right)\Picurlkdual\right)\solw_N.
\end{align*}
With this, and due to 
\begin{align*}
		\dualterm{\nabla\xi_N}{\left[\Picurlkhdual-\Picurlkdual\right]\solw_N} = 0
\end{align*}
we obtain
\begin{align*}
		\norm{\left[\Picurlkhdual-\Picurlkdual\right]\solw_N}{\knormplus}^2&\leq C\left|\dualterm{\left[\Picurlkhdual-\Picurlkdual\right]\solw_N}{\left[\Picurlkhdual-\Picurlkdual\right]\solw_N}\right|\\
																		 & = C\left|\dualterm{\left[\Picurlcom-\operatorname{I}\right]\Picurlkdual\solw_N}{\left[\Picurlkhdual-\Picurlkdual\right]\solw_N}\right| \\
																		 &\leq C\norm{\left[\Picurlcom-\operatorname{I}\right]\Picurlkdual\solw_N}{\knormplus}\norm{\left[\Picurlkhdual-\Picurlkdual\right]\solw_N}{\knormplus},
\end{align*}
where the generic constant $C>0$ depends only on $\varepsilon$ and $\zeta$. In particular, the above inequality implies
\begin{align}\label{codi2}
\norm{\left[\Picurlkhdual-\Picurlkdual\right]\solw_N}{\knormplus}\leq C\norm{\left[\Picurlcom-\operatorname{I}\right]\Picurlkdual\solw_N}{\knormplus}.
\end{align}

At this point we need some regularity properties of $\Picurlkdual\solw_N$. We notice that by construction of $\Picurlkdual$ and \eqref{dualtermdual} there holds
\begin{align*}
		\forall\psi\in\Honeimp:\ \dualtermdual{\Picurlkdual\solw_N}{\nabla\psi} = 0.
\end{align*}
That is, $\Picurlkdual\solw_N$ belongs to the auxiliary space
\begin{align*}
		\auspace:=\left\{\solw\in\HXI\ |\ \diverg\varepsilon^H\solw = 0\ {\rm in}\ \Omega, \quad {\rm and} \quad \overline{k}\varepsilon^H\solw\cdot\soln+i\diverg_{\Gamma}\zeta^H\solw_T=0\ {\rm on}\ \Gamma\right\}
\end{align*}
for which there holds the following result:
\begin{lemma}\label{auspacereg}
		Let $\Gp=\geom$ be an $\ana$-partition and suppose that the coefficient $\varepsilon\in\anatenspw{\Gp}$ satisfies \eqref{coerciveeps} and that $\zeta\in\anatens{\Gamma}$ satisfies \eqref{coercivezeta} and \eqref{zetaproperty}. Then, for all $k\in\Coone$ there holds $\auspace\subseteq\PVSHM{1}(\Gp)$ and 
		\begin{align*}
				\norm{\solv}{\PVSHM{1}(\Gp)}\leq C\left(\norm{\curl\solv}{\VSL(\Omega)}+|k|\norm{\solv}{\VSL(\Omega)}\right)
		\end{align*}
		for all $\solv\in\auspace$, where $C>0$ depends only on $\Gp$, $\varepsilon$ and $\zeta$.
\end{lemma}
\begin{fatproof}
		From the definition of $\auspace$ we have that $\solv\in\auspace$ satisfies $\diverg\varepsilon^H\solv =0$ as well as $\curl\solv\in\VSL(\Omega)$ and $\diverg_{\Gamma}\zeta^H\solv_T = -i\overline{k}\varepsilon^H\solv\cdot\soln\in\SHM{-1/2}(\Gamma)$. Furthermore, we note that $\curl_{\Gamma}\solv_T = \curl\solv\cdot\soln\in\SHM{-1/2}(\Gamma)$ and apply \cite[Lem.~3.6]{MaxwellMyself} to infer $\solv_T\in\VSHM{1/2}(\Gamma)$ with
$$ \norm{\solv_T}{\VSHM{1/2}(\Gamma)}\leq C\left(\norm{\curl\solv}{\VSL(\Omega)}+|k|\norm{\solv}{\VSL(\Omega)}\right).$$
The claim then follows from \cite[Thm.~2.7]{MaxwellMyself}.
\end{fatproof}

From \eqref{codi1} and \eqref{codi2} we infer
\begin{align*}
	|T_1|\leq\left| \dualterm{\sole_N}{(\Picurlkhdual-\Picurlkdual)\solw_N}\right|\leq C\norm{\sole_N}{\knormplus}\norm{\left[\Picurlcom-\operatorname{I}\right]\Picurlkdual\solw_N}{\knormplus},
\end{align*}
hence, together with Lemma~\ref{auspacereg}, 
\begin{align}\label{estimateT1}
		\begin{split}
				|T_1|&\leq C \norm{\sole_N}{\knormplus}\etaone\norm{\Picurlkdual\solw_N}{\PVSHM{1}(\Gp)}\\
		&\leq C\norm{\sole_N}{\knormplus}\etaone\norm{\Picurlkdual\solw_N}{\HXIK} \\
		&\leq C\norm{\sole_N}{\knormplus}\etaone\norm{\solw_N}{\HXIK}
\end{split}
\end{align}
for 
$$ \etaone := \sup_{\substack{\solw\in\auspace\setminus\{{\bf 0}\} \\ \curl\solw\in\curl\Xh}}\frac{\norm{\solw-\Picurlcom \solw}{\knormplus}}{\norm{\solw}{\PVSHM{1}(\Gp)}}.$$
It remains to estimate $\eta_1$ explicitly in $h, k$ and $p$. Before doing that we derive an analogous estimate as \eqref{estimateT1} for the term $T_2$ from \eqref{decompositionconsistency}.

\subsection{The duality argument}

We aim to provide a similar estimate as \eqref{estimateT1} for the term $T_2$ from \eqref{decompositionconsistency}. We follow \cite[Sec.~9]{MaxwellImpedanceMelenk} and for an arbitrary but fixed $\solw_N\in\Xh$ we set 
\begin{align}\label{defvzero}
\solv_0:=\Picurlkdual\solw_N.
\end{align}
We notice that due to \eqref{discretehelmholtz} and \eqref{defprop} there holds 
\begin{align*}
		\forall\psi\in\Honeimp:\ \dualtermdual{\Picurlkdual\solw_N}{\nabla\psi}^H =0,
\end{align*}
hence $\solv_0\in\auspace$. 
We pose the dual problem of finding $\solz\in\HXI$ such that
\begin{align}\label{dualproblem}
		\forall\solw\in\HXI:\ A_k(\solw,\solz) = \dualterm{\solw}{\solv_0}
\end{align}
and notice that \eqref{dualproblem} is equivalent to finding $\solz\in\HXI$ such that
\begin{align*}
		\forall\solw\in\HXI:\ A_{-\overline{k}}^H(\solz,\solw) = \dualtermdual{\solv_0}{\solw},
\end{align*}
where 
\begin{align*}
		A_{-\overline{k}}^H(\solz,\solw) := \SCP{\mu^{-H}\curl \solz}{\curl\solw}{\VSL(\Omega)}-\overline{k}^2\SCP{\varepsilon^H\solz}{\solw}{\VSL(\Omega)}+i\overline{k}\SCP{\zeta^H\solz_T}{\solw_T}{\VSL_T(\Gamma)}.
\end{align*}
Therefore, solving the dual problem \eqref{dualproblem} is equivalent to solve the primal problem
\begin{equation}\label{Maxwelldual}
\begin{alignedat}{2}
		\curl\mu^{-H}\curl\solz-\overline{k}^2\varepsilon^H\solz &= \overline{k}^2\varepsilon^H\solv_0\quad &&{\rm in}\quad \Omega, \\
		\left(\mu^{-H}\curl\solz\right)\times\soln+i\overline{k}\zeta^H\solz_T &= -i\overline{k}\zeta^H\Pi_T\solv_{0} \quad &&{\rm on}\quad  \Gamma,
\end{alignedat}
\end{equation}

We denote the solution of \eqref{Maxwelldual} by $\dualsol\solv_0$. A calculation shows that $\solv_0\in\auspace$ implies $\dualsol\solv_0\in\auspace$, hence \eqref{defvzero} and Galerkin orthogonality prove that for arbitrary $\solv_N,\solw_N\in\Xh$ we have 
\begin{align}\label{T2estimate}
		\begin{split}
		|T_2| = \left| \dualterm{\sole_N}{\Picurlkdual\solw_N}\right| = \left| \dualterm{\sole_N}{\solv_0}\right| &= \left|A_k(\sole_N, \dualsol\solv_0-\solv_N)\right|\\
																												  &\leq \etax \norm{\sole_N}{\HXIK}\norm{\solv_0}{\HXIK},
\end{split}
\end{align}
where 
\begin{align}\label{defetax}
		\etax:=\sup_{\solv_0\in\auspace\setminus\{{\bf 0}\}}\inf_{\solv_N\in\Xh}\frac{\norm{\dualsol\solv_0-\solv_N}{\HXIK}}{\norm{\solv_0}{\HXIK}}.
\end{align}

It remains to find $h,p,k$-explicit upper bounds for $\etaone$ and $\etax$.

\subsection{Wavenumber-explicit $hp$-analysis of $\etaone$ and $\etax$}

We follow the arguments from \cite[Sec.~9]{MaxwellImpedanceMelenk} to perform an $h,p,k$-explicit analysis of $\etaone$ and $\etax$. 
For $\etaone$ we may readily apply the arguments 
from \cite[Lem.~8.6]{MaxwellTransparentMelenk} and \cite[Lem.~9.4, (9.18)]{MaxwellImpedanceMelenk} and conclude
\begin{align}\label{estimateeta1}
\etaone\leq C'\left(\frac{|k|h}{p}\right)^{1/2}\left(1+\left(\frac{|k|h}{p}\right)^{1/2}\right)
\end{align}
for a constant $C'>0$ depending only on $\Gp$, $\mu^{-1}, \varepsilon, \zeta$ and the parameters $C,\gamma>0$ from Assumption~\ref{Ttrafo}. 

\bigskip

For $\etax$, we proceed as follows: According to Theorem~\ref{Mainresult1} we may write $\dualsol\solv_0 = \soluH+\soluA+\overline{k}^{-1}\nabla\psi$ with 
\begin{align*}
		\norm{\soluH}{\PVSHM{2}(\Gp)}&\leq C|k|\left(\norm{\solv_0}{\PVSHM{1}(\Gp)}+\norm{\Pi_T\solv_{0}}{\VSHM{1/2}_T(\Gamma)}\right)\leq C|k|\norm{\solv_0}{\HXIK}, \\
				\norm{\psi}{\PSHM{2}(\Gp)}&\leq C|k|\norm{\Pi_T\solv_{0}}{\VSHM{1/2}_T(\Gamma)}\leq C|k|\norm{\solv_0}{\HXIK}
\end{align*}
for a constant $C>0$ depending only on $\Gp$, $\mu^{-1}, \varepsilon$ and $\zeta$, and
\begin{align*}
		\forall\ell\in\N_0:\ \norm{\soluA}{\PVSHM{\ell}}&\leq \omega \left[|k|^2\rho(k)+|k|\right]\left(\norm{\solv_0}{\PVSHM{1}(\Gp)}+\norm{\Pi_T\solv_{0}}{\VSHM{1/2}_T(\Gamma)}\right)M^{\ell}(\ell+|k|)^{\ell} \\
								   &\leq \omega \left[|k|^2\rho(k)+|k|\right]M^{\ell}(\ell+|k|)^{\ell} \norm{\solv_0}{\HXIK},
\end{align*}
where the constants $\omega\geq0, M>0$ depend only on $\Gp$, $\mu^{-1}$, $\varepsilon$ and $\zeta$. Under the assumption \eqref{rhoalgebra} we may further estimate 
\begin{align*}
	\forall\ell\in\N_0:\ \norm{\soluA}{\PVSHM{\ell}}\leq\omega |k|^{\theta+2}M^{\ell}(\ell+|k|)^{\ell} \norm{\solv_0}{\HXIK},
\end{align*}
where $M>0$ is as before and $\omega\geq 0$ now additionally depends on the constant $C_{\setu}$ from \eqref{rhoalgebra}.

We return to \eqref{defetax} and choose $\solv_N:= \Picurlhp\soluH+\Picurlhp\soluA+\overline{k}^{-1}\nabla\Pigradhp\psi$. 
By repeating the arguments from \cite[(9.29), (9.30)]{MaxwellImpedanceMelenk} we arrive at
\begin{align}\label{approxsoluH}
		\norm{\soluH-\Picurlhp\soluH}{\HXIK}\leq C\frac{h|k|}{p}\left(1+\left(\frac{h|k|}{p}\right)^{1/2}+\frac{h|k|}{p}\right)\norm{\solv_0}{\HXIK},
\end{align}
and from Proposition~\ref{defpicurlhp} we get under the mild resolution condition
\begin{align}\label{rescondmild}
h+|k|h/p\leq C_1
\end{align}
that
\begin{align}\label{approxsoluA}
		\norm{\soluA-\Picurlhp\soluA}{\HXIK}\leq C\omega|k|^{\theta+3}\left(\left(\frac{h}{h+\sigma}\right)^p+\left(\frac{|k|h}{\sigma p}\right)^p\right)\norm{\solv_0}{\HXIK}.
\end{align}
Finally, using Proposition~\ref{defprigradhp} we argue as in \cite[(9.28)]{MaxwellImpedanceMelenk} that 
\begin{align}\label{approxpsi}
		\norm{\overline{k}^{-1}\left(\nabla\psi-\nabla\Pigradhp\psi\right)}{\HXIK}\leq C\left(\frac{|k|h}{p}\right)^{1/2}\left(1+\left(\frac{|k|h}{p}\right)^{1/2}\right)\norm{\solv_0}{\HXIK}.
\end{align}

Let us suppose that there holds \eqref{rescondmild} for some $C_1>0$ as well as $|k|h/p\leq 1$. Then, we may combine \eqref{approxsoluH}, \eqref{approxsoluA} and \eqref{approxpsi} with \eqref{defetax} to obtain
\begin{align}\label{estimate1etak}
		\etax \leq C'\left[\left(\frac{|k|h}{p}\right)^{1/2}+|k|^{\theta+3}\left(\left(\frac{h}{h+\sigma}\right)^p+\left(\frac{|k|h}{\sigma p}\right)^p\right)\right],
\end{align}
where the constants $C',\sigma>0$ depend only on $\Gp$, $\mu^{-1}, \varepsilon, \zeta$, $C_1$, the parameters $C,\gamma>0$ from Assumption~\ref{Ttrafo} and the constant $C_{\setu}$ from the assumption \eqref{rhoalgebra}.

\bigskip

From the estimate \eqref{estimateeta1} we immediately conclude that we can make $\etaone$ arbitrarily small by just choosing $h$ and $p$ such that $h|k|/p$ is sufficiently small. For $\etax$, however, the situation is more difficult as the estimate \eqref{estimate1etak} is more involved. However, based on \eqref{estimate1etak} we may follow the proof of \cite[Lem.~9.5]{MaxwellImpedanceMelenk} almost verbatim to obtain the following result:

\begin{lemma}\label{epslemma}
		Let $\Gp=\geom$ be an $\ana$-partition and let $\mu^{-1},\varepsilon\in\anatenspw{\Gp}$, $\zeta\in\anatens{\Gamma}$ and $\setu\subseteq\Coone$ be such that Assumption~\ref{assumption2} and Assumption~\ref{assumption3} hold true. Furthermore, let $\T$ be a regular and shape regular (curvilinear) triangulation of $\Omega$ as described in Section~\ref{discretizationsec} and suppose that $\T$ satisfies Assumption~\ref{Ttrafo}.

		Then, there holds the following: For every $\beta,C_2>0$ there exists a constant $C_3>0$ depending only on $\beta,C_2,\Gp,\mu^{-1},\varepsilon,\zeta$, the parameters $C,\gamma>0$ of Assumption~\ref{Ttrafo} and the constants $C_{\setu}>0$, \revision{$\theta\geq 0$} from \eqref{rhoalgebra} such that 
		\begin{align*}
				1+C_2\log |k|\leq p\quad {\rm and}\quad \frac{|k|h}{p}\leq C_3\quad {\rm imply}\quad \etax \leq \beta
		\end{align*}
		for all $k\in\setu$.
\end{lemma}

\subsection{Proof of Theorem~\ref{Mainresult2}}

At this point we have allocated all necessary tools to prove Theorem~\ref{Mainresult2}. 
Let $\solu\in\HXI$ be the solution of Maxwell's equations \eqref{Maxwellorig} and suppose that $\solu_N\in\Xh$ is a discrete vector field such that $\bfa(\solu-\solu_N, \solv_N)=0$ for all $\solv_N\in\Xh$. Our aim is to apply Proposition~\ref{quasioptimalprop}, however, in order to do that we must first verify $\delta_k(\solu-\solu_N)<\cco^{-1}$, where $\cco>0$ is the ellipticity constant from \eqref{coercivemu} and \eqref{ezcoerceps}-\eqref{ezcoerczeta}.

By definition of $\delta_k(\solu-\solu_N)$, the splitting \eqref{decompositionconsistency} and the estimates \eqref{estimateT1} and \eqref{T2estimate} we have
\begin{align*}
\delta_k(\solu-\solu_N) \leq C\left(\etaone+\etax\right),
\end{align*}
hence we may use \eqref{estimateeta1} and Lemma~\ref{epslemma} to obtain the following:

\begin{lemma}\label{deltaepsilon}
		Let $\Gp=\geom$ be an $\ana$-partition and let $\mu^{-1},\varepsilon\in\anatenspw{\Gp}$, $\zeta\in\anatens{\Gamma}$ and $\setu\subseteq\Coone$ be such that Assumption~\ref{assumption2} and Assumption~\ref{assumption3} hold true. Furthermore, let $\T$ be a regular and shape regular (curvilinear) triangulation of $\Omega$ as described in Section~\ref{discretizationsec} below and suppose that $\T$ satisfies Assumption~\ref{Ttrafo}.

Furthermore, let $\solu\in\HXI$ be the solution of Maxwell's equations \eqref{Maxwellorig} and suppose that $\solu_N\in\Xh$ is such that $\bfa(\solu-\solu_N, \solv_N) = 0$ for all $\solv_N\in\Xh$. 

Under these conditions, let $\beta,C_2>0$ be arbitrary. Then, there exists a constant $C_3>0$ depending only on $\beta, C_2, \Gp$, $\mu^{-1}, \varepsilon, \zeta$, the parameters $C,\gamma>0$ from Assumption~\ref{Ttrafo} and the constants $C_{\setu}>0$, \revision{$\theta\geq 0$} from \eqref{rhoalgebra} such that
\begin{align*}
		1+C_2\log |k|\leq p \quad{\rm and}\quad \frac{|k|h}{p}\leq C_3 \quad {\rm implies} \quad 
\delta_k(\solu-\solu_N)<\beta
\end{align*}
for all $k\in\setu$.
\end{lemma}
\begin{fatproof}
		Follows from \eqref{estimateeta1} and Lemma~\ref{epslemma}.
\end{fatproof}

We conclude this section proving Theorem~\ref{Mainresult2}.

\medskip 

\begin{fatproofmod}{Theorem~\ref{Mainresult2}}
		We first show that under the present assumptions on $h$ and $p$, the discrete problem \eqref{discreteproblem} has a unique solution $\solu_N\in\Xh$. We follow the arguments from \cite[Thm.~4.15]{MaxwellTransparentMelenk} and consider uniqueness first. Let $\solw_N\in\Xh$ be such that $\bfa(\solw_N,\solv_N) = 0$ for all $\solv_N\in\Xh$. We notice that this is the discretization of the continuous problem $\bfa(\solw,\solv)=0$ for all $\solv\in\HXI$, which according to Assumption~\ref{assumption3} has the unique solution $\solw = 0$.

		By Lemma~\ref{deltaepsilon} we may choose $C_3>0$ sufficiently small such that $2\delta_k(\solw-\solw_N)\leq \cco^{-1}$. Hence, by Proposition~\ref{quasioptimalprop} there holds
		\begin{align*}
				\norm{\solw_N}{\HXIK} = \norm{\solw-\solw_N}{\HXIK}\leq C(\cco^{-1}+2)\inf_{\solv_N\in\Xh}\norm{\solv_N}{\HXIK} = 0.
		\end{align*}
		Hence, $\solw_N=0$ and therefore solutions of \eqref{discreteproblem} are unique. Furthermore, since the discrete problem \eqref{discreteproblem} is finite-dimensional with as many equations as unknowns, uniqueness implies existence. That is, under the present assumptions on $h$ and $p$ the discrete problem \eqref{discreteproblem} has a unique solution.
		
		\medskip

		It remains to prove the quasi-optimality estimate. Let $\eta\in (0,\cco^{-1})$ and $C_2>0$ be arbitrary. Furthermore, let $\solu\in\HXI$ be the solution of Maxwell's equations \eqref{Maxwellorig} and suppose that $\solu_N\in\Xh$ is the solution of the corresponding discrete problem \eqref{discreteproblem}. According to Lemma~\ref{deltaepsilon} we may choose $C_3>0$ such that 
\begin{align*}
		1+C_2\log |k|\leq p \quad{\rm and}\quad \frac{|k|h}{p}\leq C_3 \quad {\rm implies} \quad 
\delta_k(\solu-\solu_N)<\eta
\end{align*}
for all $k\in\setu$. Hence, Proposition~\ref{quasioptimalprop} asserts
		\begin{align*}
				\norm{\solu-\solu_N}{\HXIK}\leq C\cco\frac{1+\eta}{1-\cco\eta}\ \inf_{\solw_N\in\Xh}\norm{\solu-\solw_N}{\HXIK},
		\end{align*}
		which concludes the proof.
\end{fatproofmod}

	\section{Numerical experiments}\label{numerics}

To conclude this work, we illustrate the results of Theorem~\ref{Mainresult2} by two numerical experiments. For the computations we used the software package NGSolve \cite{ngsolve}, \cite{netgen} and N\'{e}d\'elec type II elements, that is, full polynomial spaces. 
\begin{remark}
		The analysis of this work was carried out for N\'{e}d\'{e}lec type I elements, however the involved element-by-element interpolation operators $\Picurlcom$ and $\Pidivcom$ from Proposition~\ref{Ipprop} are also available for type II elements, cf. \cite[Sec.~4.8]{DissRojik}. For this reason our analysis can be carried over to N\'{e}d\'{e}lec type II elements \revtwo{almost verbatim}.
\end{remark}

\begin{figure}[h]
		\begin{minipage}{0.5\textwidth}
				\begin{tikzpicture}
  \begin{loglogaxis}[
    width=1\textwidth, height=6cm,     
    grid = major,
    grid style={dashed, gray!30},
    axis background/.style={fill=white},
    xlabel=$N_k\sim$ DOFs per wavelength,
	ylabel=rel. error in $\norm{\cdot}{\operatorname{curl},k}$,
    legend style={at={(0.02,0.02)},anchor=south west,font=\footnotesize},
    xtick = {},
    ]
      
    \addplot+[solid,mark=square,mark size=2pt,mark options={line width=1.0pt}] table    
    [
    x index = 1,
    y index = 2,
	col sep=comma,
    ]{plotdata_ball/data_p1_k10.csv};    
    \addlegendentry{$k=10$}

    \addplot+[solid,mark=star,mark size=2pt,mark options={line width=1.0pt}] table    
    [
    x index=1,
    y index = 2,
	col sep=comma,
    ]{plotdata_ball/data_p1_k20.csv};    
    \addlegendentry{$k=20$}

    \addplot+[solid,mark=triangle,mark size=2pt,mark options={line width=1.0pt}] table    
    [
    x index=1,
    y index=2,
	col sep=comma,
    ]{plotdata_ball/data_p1_k30.csv};    
    \addlegendentry{$k=30$}
    

			\addplot [black,dashed ] expression [domain=1:25, samples=15] {0.4*x^(-1)}; \addlegendentry{$\mathcal{O}(h)$} 
  \end{loglogaxis} 
\end{tikzpicture}
		\end{minipage}
		\begin{minipage}{0.5\textwidth}
				\begin{tikzpicture}
  \begin{loglogaxis}[
    width=1\textwidth, height=6cm,     
    grid = major,
    grid style={dashed, gray!30},
    axis background/.style={fill=white},
    xlabel=$N_k\sim$ DOFs per wavelength,
	ylabel=rel. error in $\norm{\cdot}{\operatorname{curl},k}$,
    legend style={at={(0.02,0.02)},anchor=south west,font=\footnotesize},
    xtick = {},
    ]
      
    \addplot+[solid,mark=square,mark size=2pt,mark options={line width=1.0pt}] table    
    [
    x index = 1,
    y index = 2,
	col sep=comma,
    ]{plotdata_ball/data_p2_k10.csv};    
    \addlegendentry{$k=10$}

    \addplot+[solid,mark=star,mark size=2pt,mark options={line width=1.0pt}] table    
    [
    x index=1,
    y index = 2,
	col sep=comma,
    ]{plotdata_ball/data_p2_k20.csv};    
    \addlegendentry{$k=20$}

    \addplot+[solid,mark=triangle,mark size=2pt,mark options={line width=1.0pt}] table    
    [
    x index=1,
    y index=2,
	col sep=comma,
    ]{plotdata_ball/data_p2_k30.csv};    
    \addlegendentry{$k=30$}
    

			\addplot [black,dashed ] expression [domain=1:17, samples=15] {0.4*x^(-2)}; \addlegendentry{$\mathcal{O}(h^2)$} 
  \end{loglogaxis} 
\end{tikzpicture}
		\end{minipage}
		\caption{\revtwo{FEM convergence for} piecewise constant coefficients, left: $p=1$, right: $p=2$.}
		\label{fig:ex1}
\end{figure}
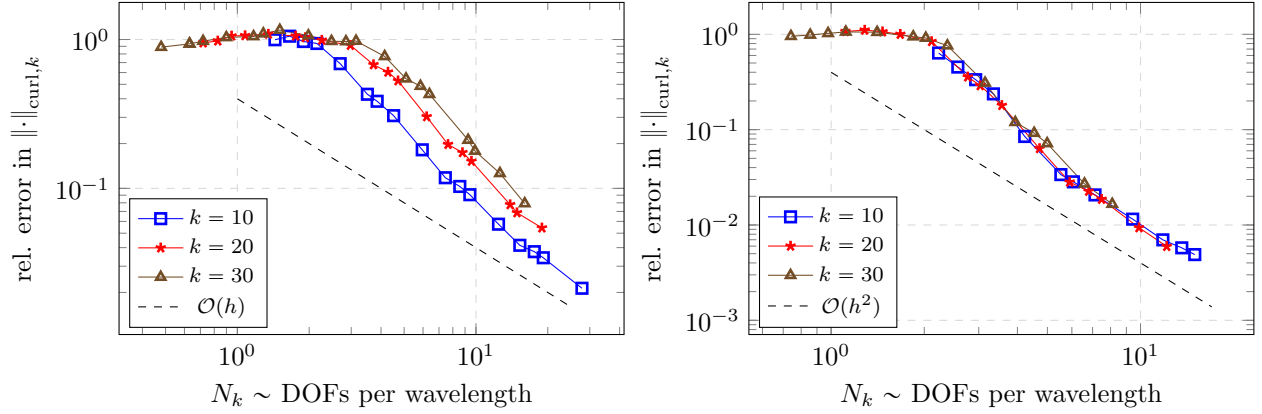

Our first example is set on an $\ana$-partition $\Gp=\{\Omega, \Gp_1,\Gp_2\}$ with 
$$\Omega=\mathcal{B}_1(0), \quad \Gp_1 = \mathcal{B}_{1/2}(0)\quad\text{and}\quad \Gp_2 = \mathcal{B}_1(0)\setminus\Gp_1.$$
We prescribe piecewise constant coefficients
\begin{align*}
		\mu^{-1} = \begin{pmatrix}
				3 & 1 & 0 \\
				1 & 3 & 1 \\
				0 & 1 & 3
		\end{pmatrix}
		\quad\text{and}\quad 
\varepsilon = \begin{pmatrix}
				2 & 1 & 0 \\
				1 & 2 & 0 \\
				0 & 0 & 3
		\end{pmatrix}
		\quad\text{on}\quad \Gp_2,
\end{align*}
and $\mu^{-1}=\varepsilon = 1$ on $\Gp_1$. With these choices of $\Gp$, $\mu^{-1}$ and $\varepsilon$ and with $\solf(x,y,z):=(z,0,0)^T$ we consider the equations 
\begin{equation}\label{numerics:equation}
\begin{alignedat}{2}			\curl\mu^{-1}\curl\solu-k^2\varepsilon\solu &= \solf \quad &&\text{in}\ \Omega, \\
		\mu^{-1}\curl\solu\times\soln-ik\solu_T &=0 \quad &&\text{on}\ \Gamma.
	\end{alignedat}
\end{equation}
For the values $k\in\{10,20,30\}$ we solve these equations using $hp$-FEM with fixed polynomial degrees $p\in\{1,2\}$. Furthermore, for every value $k\in\{10,20,30\}$ we compute a reference solution using polynomial order $p=4$.
In Figure~\ref{fig:ex1} we plot the resulting relative errors in the norm 
$$\norm{\cdot}{\operatorname{curl},k}^2:=\norm{\curl\cdot}{\VSL(\Omega)}^2+|k|^2\norm{\cdot}{\VSL(\Omega)}^2$$ 
over the quantity
$$N_k:=\frac{\operatorname{DOF}^{1/3}}{|k|}$$
where $\operatorname{DOF}$ is the dimension of the discrete space $\Xh$. We remark that the quantity $N_k$ is a measure for the \revtwo{number} of degrees of freedom per wavelength; indeed, in the case of homogeneous coefficients $\mu^{-1}=\varepsilon=1$ the \revtwo{number of} degrees of freedom per wavelength \revtwo{is} given by
$$\text{DOFs per wavelength}=\frac{2\pi \operatorname{DOF}^{1/3}}{|\Omega|^{1/3}|k|}.$$
In the presence of heterogeneous coefficients the DOFs per wavelength are not clearly defined, however $N_k$ is still an estimate of them.

We observe the expected $\mathcal{O}(h^{p})$ behavior, and we further observe that for $p=2$ the finite element error reaches its asymptotic behavior for lower values of $N_k$. This indicates that for the choice $p=2$ the $hp$-FEM suffers less from pollution compared to $p=1$.

\medskip

\begin{figure}[h]
		\begin{minipage}{0.5\textwidth}
				\begin{tikzpicture}
  \begin{loglogaxis}[
    width=1\textwidth, height=6cm,     
    grid = major,
    grid style={dashed, gray!30},
    axis background/.style={fill=white},
    xlabel=$N_k\sim$ DOFs per wavelength,
	ylabel=rel. error in $\norm{\cdot}{\operatorname{curl},k}$,
    legend style={at={(0.02,0.02)},anchor=south west,font=\footnotesize},
    xtick = {},
    ]
      
    \addplot+[solid,mark=square,mark size=2pt,mark options={line width=1.0pt}] table    
    [
    x index = 1,
    y index = 2,
	col sep=comma,
    ]{plotdata_brick/data_p1_k10.csv};    
    \addlegendentry{$k=10$}

    \addplot+[solid,mark=star,mark size=2pt,mark options={line width=1.0pt}] table    
    [
    x index=1,
    y index = 2,
	col sep=comma,
    ]{plotdata_brick/data_p1_k20.csv};    
    \addlegendentry{$k=20$}

%
	\addplot+[solid,mark=triangle,mark size=2pt,mark options={line width=1.0pt}] table    
    [
    x index=1,
    y index=2,
	col sep=comma,
    ]{plotdata_brick/data_p1_k40.csv};    
    \addlegendentry{$k=40$}


			\addplot [black,dashed ] expression [domain=1:20, samples=15] {0.2*x^(-1)}; \addlegendentry{$\mathcal{O}(h)$} 
  \end{loglogaxis} 
\end{tikzpicture}
		\end{minipage}
		\begin{minipage}{0.5\textwidth}
				\begin{tikzpicture}
  \begin{loglogaxis}[
    width=1\textwidth, height=6cm,     
    grid = major,
    grid style={dashed, gray!30},
    axis background/.style={fill=white},
    xlabel=$N_k\sim$ DOFs per wavelength,
	ylabel=rel. error in $\norm{\cdot}{\operatorname{curl},k}$,
    legend style={at={(0.02,0.02)},anchor=south west,font=\footnotesize},
    xtick = {},
    ]
      
    \addplot+[solid,mark=square,mark size=2pt,mark options={line width=1.0pt}] table    
    [
    x index = 1,
    y index = 2,
	col sep=comma,
    ]{plotdata_brick/data_p2_k10.csv};    
    \addlegendentry{$k=10$}

    \addplot+[solid,mark=star,mark size=2pt,mark options={line width=1.0pt}] table    
    [
    x index=1,
    y index = 2,
	col sep=comma,
    ]{plotdata_brick/data_p2_k20.csv};    
    \addlegendentry{$k=20$}

%
	\addplot+[solid,mark=triangle,mark size=2pt,mark options={line width=1.0pt}] table    
    [
    x index=1,
    y index=2,
	col sep=comma,
    ]{plotdata_brick/data_p2_k40.csv};    
    \addlegendentry{$k=40$}


			\addplot [black,dashed ] expression [domain=1:20, samples=15] {0.1*x^(-2)}; \addlegendentry{$\mathcal{O}(h^2)$} 
  \end{loglogaxis} 
\end{tikzpicture}
		\end{minipage}
		\caption{\revtwo{FEM convergence for} smooth coefficients, left: $p=1$, right: $p=2$.}
		\label{fig:ex2}
\end{figure}
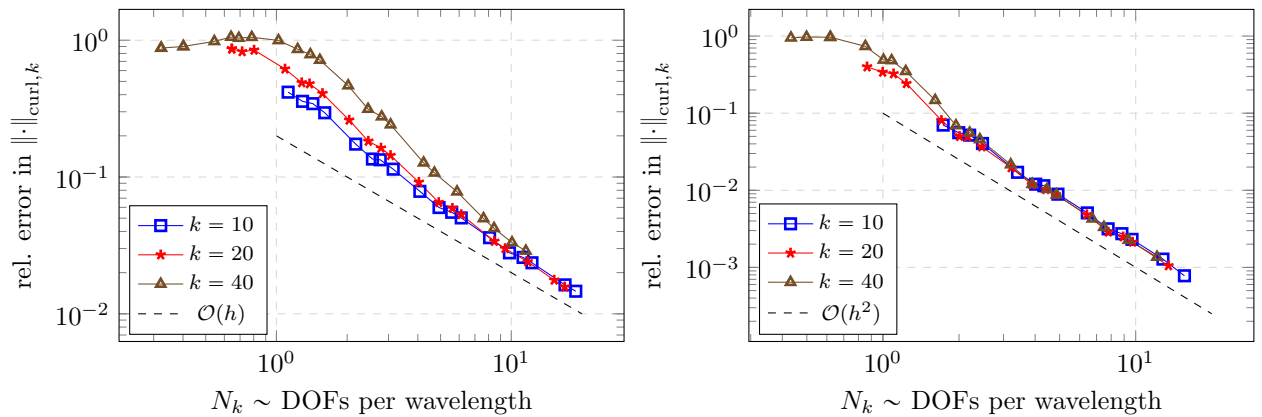

In our second example we consider the unit cube $\Omega= (0,1)^3$ and solve \eqref{numerics:equation} for
\begin{align*}
		\mu^{-1}(x,y,z) = 1+x^2,\quad \varepsilon=1 \quad\text{and}\quad \solf(x,y,z) = (z,2x,0)^T\quad\text{in}\ \Omega.
\end{align*}
For $k\in\{10,20,40\}$ we compute reference solutions using polynomial degree $p=4$, and using these reference solutions we can calculate the finite element error of $hp$-FEM for $p=1$ and $p=2$. In Figure~\ref{fig:ex2} we plot the relative error measured in the $\norm{\cdot}{\operatorname{curl},k}$ over $N_k$, and we again observe the expected $\mathcal{O}(h^p)$ convergence. As in the first example we observe that for $p=2$ the finite element error reaches its asymptotic behaviour for lower values of $N_k$ which indicates that the pollution effect is mitigated. Let us remark that this example is actually not covered by Theorem~\ref{Mainresult2} since the boundary of $\Omega$ is only piecewise analytic, however the numerical evidence suggests that Theorem~\ref{Mainresult2} is valid nevertheless.

\newpage

\section*{Acknowledgement}
    JMM acknowledges funding by the Austrian Science Fund (FWF) under
    grant F65 ``taming complexity in partial differential systems'' (\href{https://doi.org/10.55776/F65}{DOI:10.55776/F65}).

	
	\bibliographystyle{plain}
	\bibliography{bibliog}
\end{document}